\documentclass[12pt]{article}
\usepackage{latexsym, amssymb, amsmath, amscd, amsfonts, epsfig, graphicx, colordvi,verbatim,ifpdf}
\usepackage{amsfonts, amsmath, amssymb}
\usepackage{amssymb,amsfonts,amsmath,latexsym,epsfig,cite, psfrag,eepic,color}
\usepackage{amscd,graphics}
\usepackage{latexsym, amssymb,  amsmath,amscd, amsfonts, epsfig, graphicx, colordvi,amsthm}

\usepackage{graphicx}
\usepackage{color}
\usepackage{ifpdf}
\usepackage{fancybox}

\usepackage{float}

\newtheorem{thm}{Theorem}[section]
\newtheorem{pro}[thm]{Proposition}
\newtheorem{conj}[thm]{Conjecture}

\newtheorem{lem}[thm]{Lemma}

\def\pf{\noindent{\it Proof.} }
\def\qed{\nopagebreak\hfill{\rule{4pt}{7pt}}
\medbreak}

\numberwithin{equation}{section}

\def\qed{\nopagebreak\hfill{\rule{4pt}{7pt}}
\medbreak}

\newlength{\boxedparwidth}
  {\begin{center} \begin{tabular}{|@{\hspace{.315in}}c@{\hspace{.15in}}|}
                  \hline \\ \begin{minipage}[t]{\boxedparwidth}
                  \setlength{\parindent}{.25in}}%
  {\end{minipage} \\ \\ \hline \end{tabular} \end{center}}

\begin{document}
\begin{center}

 {\Large \bf Proof of the Andrews-Dyson-Rhoades Conjecture

  on the spt-Crank}
\end{center}

\begin{center}
{  William Y.C. Chen}$^{1}$, {Kathy Q. Ji}$^{2}$ and
  {Wenston J.T. Zang}$^{3}$ \vskip 2mm

   $^{1,2,3}$Center for Combinatorics, LPMC-TJKLC\\
   Nankai University, Tianjin 300071, P.R. China\\[6pt]
   $^{1}$Center for Applied Mathematics\\
Tianjin University,  Tianjin 300072, P. R. China\\[6pt]

   \vskip 2mm

    $^1$chen@nankai.edu.cn, $^2$ji@nankai.edu.cn, $^3$wenston@mail.nankai.edu.cn
\end{center}

\vskip 6mm \noindent {\bf Abstract.} The notion of the spt-crank of a vector
partition, or an $S$-partition, was introduced by Andrews, Garvan and Liang.
 Let $N_S(m,n)$ denote the number of $S$-partitions of $n$ with spt-crank $m$. Andrews, Dyson and Rhoades conjectured that $\{N_S(m,n)\}_m$ is unimodal for any $n$, and they  showed that this conjecture is equivalent to an inequality between  the rank and the crank of ordinary partitions.   They obtained an asymptotic formula for the
   difference between   the rank and the crank of ordinary partitions,
   which implies $N_S(m,n)\geq N_S(m+1,n)$ for  sufficiently large $n$ and fixed $m$. In this paper, we introduce
   a representation of an ordinary  partition,  called
    the $m$-Durfee rectangle symbol, which is a rectangular
    generalization of the Durfee symbol introduced by Andrews.
     We give a proof of the conjecture of
      Andrews, Dyson and Rhoades by considering two cases.
      For $m\geq 1$, we construct an injection
    from the set of ordinary partitions of $n$ such that $m$ appears in the rank-set to the set of ordinary partitions of $n$ with rank  not less than $-m$.
    The case for $m=0$ requires five more injections.  We also show that this conjecture implies an inequality
between the positive rank  and  crank moments obtained by
 Andrews, Chan and Kim.

\noindent {\bf Keywords}: Inequality, rank, crank, spt-crank, Andrews' spt-function, rank moments, crank moments.

\noindent {\bf AMS Classifications}: 05A17, 11P82, 11P83.

\section{Introduction}

   In this paper, we give a proof of a conjecture of Andrews, Dyson and
   Rhoades on the spt-crank of a vector partition or an $S$-partition.
   The notion of the spt-function, called the smallest part function,
   was introduced by Andrews \cite{Andrews-2008}.  More precisely, we use
   $spt(n)$ to denote the total number of smallest parts in all partitions of $n$.
   For example, we have $spt(3)=5$, $spt(4)=10$ and $spt(5)=14$. The smallest part function  possesses many
arithmetic properties analogous to the ordinary partition function, see, for example, \cite{Andrews-2008,   Folsom-Ono-2008, Garvan, Ono-2011}.

     Andrews \cite{Andrews-2008} showed that the spt-function satisfies the following  Ramanujan type congruences:
    \begin{eqnarray}
 spt(5n+4)&\equiv& 0 \pmod 5, \label{con-5}\\[3pt]
spt(7n+5) &\equiv& 0 \pmod 7, \label{con-7}\\[3pt]
spt(13n+6) &\equiv& 0 \pmod {13}.
\end{eqnarray}

 To give combinatorial interpretations of the above congruences,
 Andrews, Garvan and Liang \cite{Andrews-Garvan-Liang-2011} introduced the   spt-crank of an $S$-partition.  Let $\mathcal{D}$ denote  the set of partitions into distinct parts and $\mathcal{P}$ denote the set of partitions. For   $\pi \in \mathcal{P}$, we use $s(\pi)$ to denote  the smallest part  of $\pi$ with the convention that $s(\emptyset)=+\infty$. Let $\ell(\pi)$ denote the number of parts of $\pi$.   Define
\[S=\{(\pi_1, \pi_2,\pi_3)\in \mathcal{D} \times \mathcal{P} \times \mathcal{P} \colon \pi_1 \neq \emptyset \text{ and }s(\pi_1)\leq \min\{s(\pi_2),s(\pi_3)\}\}.\]

For $\pi=(\pi_1,\,\pi_2,\,\pi_3) \in S$, define $|\pi|=|\pi_1|+|\pi_2|+|\pi_3|.$
Andrews, Garvan and Liang \cite{Andrews-Garvan-Liang-2011} called $\pi$ to be a $S$-partitions of $n$ if  $|\pi|=n$. The spt-crank of $\pi$, denoted  $r(\pi)$, is defined to be the difference between the number of parts of $\pi_2$ and   $\pi_3$, that is, \[ r(\pi)=\ell(\pi_2)-\ell(\pi_3).\]
For a vector partition $\pi$, we associate it with a sign
 $\omega(\pi)=(-1)^{\ell(\pi_1)-1}$. Let $N_S(m,n)$ denote the
net number of vector partitions of $n$ in $S$ with spt-crank $m$,
that is,
\begin{equation}
N_S(m,n)=\sum_{\stackrel{|\pi|=n}{r(\pi)=m}}\omega(\pi)
\end{equation}
and
\[N_S(m,t,n)=\sum_{k \equiv m \pmod{t}}N_S(k,n).\]

Andrews, Garvan and Liang \cite{Andrews-Garvan-Liang-2011}  established the following relations:
\begin{eqnarray*}
N_S(k,5,5n+4)&=& \frac{spt(5n+4)}{5}, \quad \text{for} \quad 0\leq k\leq 4, \\[3pt]
N_S(k,7,7n+5)&=& \frac{spt(7n+5)}{7}, \quad \text{for} \quad 0\leq k\leq 6,
\end{eqnarray*}
which imply the   spt-congruences \eqref{con-5} and \eqref{con-7}  respectively.

The following conjecture was posed by Andrews, Dyson and Rhoades \cite{Andrews-Dyson-Rhoades-2013}.

  \begin{conj}\label{conj-o} For $m\geq 0$ and $n\geq 0$, we have
 \[N_S(m,n)\geq N_S(m+1,n).\]
 \end{conj}

   Andrews, Dyson and Rhoades \cite{Andrews-Dyson-Rhoades-2013} showed that this conjecture is equivalent to an inequality between  the rank and the crank of ordinary partitions. Recall that the rank of an ordinary partition was
    introduced by Dyson \cite{Dyson-1944} as the largest part of the partition minus the number of parts.  The crank of an ordinary partition was defined by Andrews and Garvan \cite{Andrews-Garvan-1988} as the largest part if the partition contains no ones, and otherwise as the number of parts larger than the number of ones minus the number of ones.

 Andrews, Dyson and Rhoades \cite{Andrews-Dyson-Rhoades-2013}  found the following connection between the inequality on $N_S(m,n)$ and the inequality on the rank and the crank.

 \begin{thm}\label{relation-spt-crank}  Let $N(m,n)$ denote the number of  partitions of $n$ with rank $m$ and $M(m,n)$ denote the number of partitions of $n$ with crank $m$. Set
\[M(0,1)=-1,\quad M(-1,1)=M(1,1)=1,\quad    M(m,1)=0,\]
and define
\begin{eqnarray}
N_{\leq m}(n)&=&\sum_{|r|\leq m}N(r,n), \\[3pt]
M_{\leq m}(n)&=&\sum_{|r|\leq m}M(r,n).
\end{eqnarray}
Then for  $m\geq 0$  and $n>1$, we have
 \begin{equation}\label{equiv-e}
 N_S(m,n)- N_S(m+1,n)=\frac{1}{2}\left(N_{\leq m}(n)- M_{\leq m}(n)\right).
 \end{equation}

 \end{thm}

It is clear from (\ref{equiv-e}) that Conjecture  \ref{conj-o} is equivalent to the following  conjecture.

\begin{conj}  \label{conj}
For $m\geq0$ and $n\geq0$, we have
\begin{equation}\label{conj-e}
N_{\leq m}(n) \geq M_{\leq m}(n).
\end{equation}
\end{conj}

When $m=0$, inequality (\ref{conj-e}) was  conjectured by Kaavya \cite{Kaavya-2011}. Andrews, Dyson and Rhoades  \cite{Andrews-Dyson-Rhoades-2013} obtained the following asymptotic formula for $N_{\leq m}(n)-M_{\leq m}(n)$, which implies that  Conjecture \ref{conj} holds for fixed $m$ and sufficiently large $n$.

\begin{thm} For each $m\geq 0$, we have
\begin{equation}
(N_{\leq m}(n)-M_{\leq m}(n))\sim \frac{(2m+1)\pi^2}{192\sqrt{3}n^2}\exp\left(\pi\sqrt{\frac{2n}{3}}\right) \quad \text{as} \quad n\rightarrow \infty.
\end{equation}
\end{thm}

The main objective of this paper is to  give a proof of
Conjecture \ref{conj}. It is easy to check that Conjecture \ref{conj} holds for $n=0$ and $n=1$. To show Conjecture \ref{conj} holds for $n>1$,   we first give
 a reformulation of Conjecture \ref{conj} in terms of the rank-set. We then give an injective proof of the equivalent inequality.   Let $\lambda=(\lambda_1,\lambda_2,\ldots, \lambda_\ell)$ be an ordinary partition. Recall that the rank-set  of $\lambda$ introduced by Dyson \cite{Dyson-1989}
 is an infinite sequence
 \[ [- \lambda_1, 1-\lambda_2, \ldots, j-\lambda_{j+1},\ldots,\ell-1-\lambda_{\ell}, \ell, \ell+1,\ldots] .\]
 For example, the rank-set of $\lambda=(5,5,4,3,1)$ is $ [-5,-4,-2,0,3,5,6,7,8,\ldots].$

Dyson \cite{Dyson-1989}  also introduced the
 number $q(m,n)$  of partitions $\lambda$ of $n$ such that $m$
  appears in the rank-set of $\lambda$. For example, there are
 three partitions of $4$ whose  rank-set contains the element $1$:
 \[(4),\quad (2,1,1), \quad (1,1,1,1).\]
So we have $q(1,4)=3.$

 Dyson \cite{Dyson-1989} established the connection between
the number $q(m,n)$ and the number of partitions of $n$ with
  a bounded crank. To be more specific, let $M(\leq m,n)$ denote the number of partitions of $n$ with crank not greater than $m$.
Dyson \cite{Dyson-1989} obtained the following relation  for $n>1$,
\begin{equation}\label{dyson-eq}
M(\leq m,n)=q(m,n),
\end{equation}
 see also Berkovich and Garvan \cite[(3.5)]{Berkovich-Garvan-2002}.
Moreover, Dyson \cite{Dyson-1969, Dyson-1989} proved the following symmetries of  $N(m,n)$ and $M(m,n)$:
\begin{eqnarray}\label{sym-1}
N(m,n)&=&N(-m,n),\\[3pt]
 M(m,n)&=&M(-m,n).\label{sym-2}
\end{eqnarray}

Using relations  \eqref{dyson-eq}, \eqref{sym-1} and \eqref{sym-2}, we found the following connection between $N_{\leq m}(n)-M_{\leq m}(n)$ and $p(-m,n)-q(m,n)$, where  $p(-m,n)$ is the number of partitions of $n$ with rank not less than $-m$:

\begin{thm}\label{equi-con} For $m\geq 0$ and $n>1$, we have
\begin{equation}\label{eqn-2}
N_{\leq m}(n)-M_{\leq m}(n)=2(p(-m,n)-q(m,n)).
\end{equation}
\end{thm}

It is clear from  \eqref{eqn-2} that Conjecture \ref{conj} is equivalent to the following assertion.

\begin{thm}\label{main-equi} For $m\geq 0$ and $n\geq 1$, we have
\begin{align}\label{main-e}
q(m,n)\leq p(-m,n).
\end{align}
\end{thm}

To prove the above theorem, we first introduce a  representation of an ordinary partition, called the $m$-Durfee rectangle symbol, which is a generalization of the Durfee symbol introduced by Andrews \cite[p.48]{Andrews-2007}. Using this representation, we give characterizations  of  partitions counted by $q(m,n)$ and $p(-m,n)$. We then construct an injection from the set of partitions of $n$ such that $m$ appears in the rank-set to the set of partitions of $n$ with rank not less than $-m$.

We also find that Conjecture \ref{conj} leads to the following  inequality between the positive rank moments $\overline{N}_k(n)$  and the positive crank moments  $\overline{M}_k(n)$ obtained by Andrews, Chan and Kim \cite{Andrews-Chan-Kim},
where
  \begin{align}
\overline{N}_k(n)&=\sum_{m=1}^{+\infty}m^kN(m,n),\\[3pt]
\overline{M}_k(n)&=\sum_{m=1}^{+\infty}m^kM(m,n).
  \end{align}

\begin{thm}{\rm(\!\!\cite[Theorem 1.5]{Andrews-Chan-Kim})}\label{main-c} For  $k\geq 1$ and $n\geq 1$,  we have
 \begin{equation}\label{in-1}
 \overline{M}_k(n)>\overline{N}_k(n).
 \end{equation}
\end{thm}

Bringmann and Mahlburg \cite{Bringmann-Mahlburg-2012}
proved that the above inequality \eqref{in-1} holds
for each fixed $k\geq 1$ and sufficiently large $n$ by deriving the following asymptotic formula for $ \overline{M}_k(n)-\overline{N}_k(n)$.

  \begin{thm}For $k\geq 1$, we have
  \begin{equation}
   \overline{M}_k(n)-\overline{N}_k(n)\sim k!\zeta(k-2)(1-2^{3-k})\frac{6^{\frac{k-1}{2}}}{4\sqrt{3}\pi^{k-1}}n^{\frac{k}{2}
   -\frac{3}{2}}\exp\left(\pi \sqrt{\frac{2n}{3}}\right) \quad \text{as} \quad n\rightarrow \infty.
  \end{equation}
  Here $\zeta(s)$ denotes the Riemann $\zeta$-function.
  \end{thm}

  When $k$ is even,  inequality \eqref{in-1}   is equivalent to an  inequality of Garvan  on the ordinary rank moments $N_k(n)$ and the ordinary crank moments $M_k(n)$ introduced by  Atkin and Garvan \cite{Atkin-Garvan-03}.   For $k\geq 1$ and $n\geq 1$, Garvan \cite{Garvan-2011}  proved that
\begin{equation}\label{G-Ine}
M_{2k}(n)>N_{2k}(n),
\end{equation}
where
\begin{eqnarray*}
N_k(n)&=&\sum_{m=-\infty}^{+\infty}m^kN(m,n),\\[3pt]
M_k(n)&=&\sum_{m=-\infty}^{+\infty}m^kM(m,n).
\end{eqnarray*}

  This paper is organized as follows. In Section 2, we give a proof of Theorem \ref{equi-con}.  By Theorem \ref{equi-con}, we see that   Conjecture \ref{conj} is equivalent  to Theorem \ref{main-equi}.   In Section 3, we give the definition of a   representation of an ordinary partition, called the $m$-Durfee rectangle  symbol, which is a generalization of the Durfee symbol introduced by Andrews.   Using this new symbol,  we give characterizations  of  partitions counted by $q(m,n)$ and $p(-m,n)$.   In Section 4, we present an injective proof of Theorem \ref{main-equi} for the case $m\geq 1$. To this end, we shall build an injection from the set of partitions counted by $q(m,n)$ to the set of partitions counted by $p(-m,n)$. We divide the
  set   of partitions counted by $q(m,n)$ into six disjoint subsets $Q_i(m,n)$ ($1\leq i\leq 6$) and divide the set of partitions counted by $p(-m,n)$ into eight disjoint subsets $P_{i}(-m,n)$  ($1\leq i\leq 8$). The  injection
  consists of six injections  $\phi_i$  from the set $Q_i(m,n)$  to the set $P_{i}(-m,n)$, where $1\leq i\leq 6$.
   In Section 5, we provide a   proof of  Theorem \ref{main-equi} for the case $m=0$.   It turns out that the case $m=0$ is not  simpler than the general case $m\geq 1$.  The injections $\phi_1, \phi_2,\phi_3, \phi_4$
in Section 4 also apply to the  sets $Q_i(0,n)$, where $1\leq i\leq 4$.
We further divide $Q_5(0,n)\cup Q_6(0,n)$ into five disjoint subsets $\bar{Q}_i(0,n)$ ($1\leq i\leq 5$) and divide  $P_5(0,n)\cup P_6(0,n)$ into three disjoint subsets $\bar{P}_i(0,n)$ ($1\leq i\leq 3$).
 In addition to the two injections $\phi_5$ and $\phi_6$,
  we still need three more injections.  In Section 6, we
     demonstrate that  Theorem \ref{main-c} of Andrews, Chan and Kim can be deduced from
       Conjecture \ref{conj}.

\section{Proof of Theorem \ref{equi-con}}

In this section, we   give a proof of Theorem \ref{equi-con} which implies that
 Conjecture \ref{conj} is equivalent to Theorem \ref{main-equi}.

{\noindent \it Proof of Theorem \ref{equi-con}.} Since
\[
N_{\le m}(n)=\sum_{r=-m}^m N(r,n)
\]
and
\[
p(-m,n)=\sum_{r=-m}^{+\infty} N(r,n),
\]
we get
\[N_{\le m}(n)=p(-m,n)- \sum_{r=-\infty}^{+\infty} N(r,n)+\sum_{r=-\infty}^m N(r,n) .\]
But
 \begin{equation}\label{sum}
 \sum_{r=-\infty}^\infty N(r,n)=p(n),
 \end{equation}
so we have
 \begin{equation}\label{temp}
 N_{\le m}(n)=p(-m,n)-p(n)
 +\sum_{r=-\infty}^m N(r,n).
 \end{equation}
Replacing $r$ by $-r$ in the summation on the right-hand side of \eqref{temp}, and using the symmetry $N(m,n)=N(-m,n)$ in \eqref{sym-1}, we arrive at
\begin{eqnarray}\label{temp-2}
\sum_{r=-\infty}^m N(r,n)=\sum_{r=-m}^{+\infty} N(-r,n)=\sum_{r=-m}^{+\infty} N(r,n)=p(-m,n).
\end{eqnarray}
Substituting \eqref{temp-2} into \eqref{temp}, we obtain
\begin{equation}
N_{\le m}(n)=2p(-m,n)-p(n).\label{Nleq}
\end{equation}
Similarly, for $n>1$  we get
\begin{equation}\label{Mleq}
M_{\le m}(n)=2q(m,n)-p(n).
\end{equation}
 Subtracting \eqref{Mleq} from \eqref{Nleq} gives
\[
N_{\leq m}(n)-M_{\leq m}(n)=2(p(-m,n)-q(m,n)).
\]
This completes the proof. \qed

\section{The $m$-Durfee rectangle symbol}

In this section, we give a representation of an ordinary partition, called the
$m$-Durfee rectangle symbol, which is a generalization of the Durfee symbol introduced by Andrews \cite{Andrews-2007}.  Using this new symbol,  we give characterizations  of  partitions counted by $q(m,n)$ and $p(-m,n)$.

 Recall that the $m$-Durfee rectangle of a partition $\lambda$ introduced by Gordon and Houten \cite{Gordon-Houten-1968} is defined to be the largest $(m+j)\times j$ rectangle  contained in the Ferrers diagram of $\lambda$. An $m$-Durfee rectangle reduces to a Durfee square when $m=0$.

To give a combinatorial interpretation of the ordinary rank moments $N_{2k}(n)$,
for a partition $\lambda$ of $n$, Andrews \cite{Andrews-2007}  defined the Durfee symbol of $\lambda$ as follows
\begin{equation}
(\alpha,\beta)_D=\left(\begin{array}{cccc}
\alpha_1,&\alpha_2,&\ldots,&\alpha_s\\[3pt]
\beta_1,&\beta_2,&\ldots,&\beta_t
\end{array} \right)_{D},
\end{equation}
where $D$ is the size of the Durfee square of the Ferrers diagram of $\lambda$ and
  $\alpha$ consists of  columns to the right of the Durfee square and $\beta$ consists of  rows below the Durfee square, see Figure 3.1.
   It is clear that
  $D\geq \alpha_1\geq \alpha_2\geq \cdots \geq \alpha_s$,
$D\geq \beta_1 \geq \beta_2 \geq \cdots \geq \beta_t$, and
\[n=\sum_{i=1}^s \alpha_i+\sum_{i=1}^t \beta_i+D^2.\]
The advantage of arranging the two partitions $\alpha$ and $\beta$
in two lines lies in the fact that the rank of $\lambda$
equals the difference of the lengths of $\alpha$ and $\beta$.

  \input{Fig2.TpX}

For the partition $\lambda=(7,7,6,4,3,3,2,2,2)$ in Figure 3.1,
  the Durfee symbol of $\lambda$ is
\[(\alpha,\beta)_D= \left(\begin{array}{lllll}
3,&3,&2\\[3pt]
3,&3,&2,&2,&2
\end{array}
\right)_{4}.\]
It is easy to see that the rank of $\lambda$ is equal to $\ell(\alpha)-\ell(\beta)=-2.$

As a generalization of the Durfee symbol, the $m$-Durfee rectangle symbol
will be used in the proof of Theorem \ref{main-equi}. An $m$-Durfee rectangle symbol  is defined as follows
\begin{equation}
(\alpha,\beta)_{(m+j)\times j}=\left(\begin{array}{cccc}
\alpha_1,&\alpha_2,&\ldots,&\alpha_s\\[3pt]
\beta_1,&\beta_2,&\ldots,&\beta_t
\end{array} \right)_{(m+j)\times j},
\end{equation}
where $(m+j)\times j$ is the $m$-Durfee rectangle of the Ferrers diagram of $\lambda$ and $\alpha$ consists of  columns to the right of the $m$-Durfee rectangle and   $\beta$ consists of  rows below the $m$-Durfee rectangle, see Figure 3.2. Clearly, we have
 $m+j\geq \alpha_1\geq \alpha_2\geq \cdots \geq \alpha_s$,
$j\geq \beta_1 \geq \beta_2 \geq \cdots \geq \beta_t$ and
\[n=\sum_{i=1}^s \alpha_i+\sum_{i=1}^t \beta_i+j(m+j).\]
 \input{Fig1.TpX}
The $2$-Durfee rectangle symbol of $\lambda=(7,7,6,4,3,3,2,2,2)$ in Figure 3.2 is
\[ \left(\begin{array}{cccc }
4,&3,&3, &2\\[3pt]
3,&2,&2, &2
\end{array}
\right)_{5\times 3}.\]

Notice that for  a partition $\lambda$ with  $\ell(\lambda)\leq m$,
 there is no $m$-Durfee rectangle. In this case,
 we adopt a convention that the  $m$-Durfee rectangle has
 no columns, that is, we set $j=0$. So  $m$-Durfee rectangle symbol is  $(\lambda',\emptyset)_{m\times  0},$ where $ \lambda'$ is
  the conjugate of $\lambda$. For example, the $3$-Durfee rectangle symbol of $\lambda=(5,5,1)$ is
\[ \left(\begin{array}{ccccc }
3,&2,&2, &2,&2\\[3pt]
&
\end{array}
\right)_{3\times 0}.\]

The following
property  will be used in the next section
 to describe partitions counted by $q(m,n)$
and $p(-m,n)$ in terms of $m$-Durfee rectangle symbols.

\begin{pro}\label{pro-dufrec}Let $\lambda$ be an ordinary partition and $(\alpha,\beta)_{(m+j)\times j}$ be the $m$-Durfee rectangle symbol of $\lambda$. If $m$ appears in the rank-set of $\lambda$, then either $j=0$ or $j\geq 1$ and
\begin{equation} \beta_1=j.
\end{equation}
If the rank of $\lambda$ is not less than $-m$,   then either $j=0$ or  $j\geq 1$ and
 \begin{equation}
\ell(\beta)\leq \ell(\alpha).
 \end{equation}
  \end{pro}

\section{Proof of Theorem  \ref{main-equi} for $m\geq 1$}

 Let $Q(m,n)$ denote the set of partitions $\lambda$ of $n$ such that $m$ appears in the rank-set of $\lambda$ and   $P(-m,n)$ denote the set of partitions of $n$ with rank  not less than $-m$. Theorem \ref{main-equi} is equivalent to the following combinatorial statement.

\begin{thm}\label{thm-equi-co} For $m \geq 0$, there is an injection $\Phi$ from the set $Q(m,n)$ to the set $P(m,n)$.
\end{thm}
In this section, we give a proof of Theorem \ref{thm-equi-co} for $m\geq 1$.
The case for $m=0$ will be dealt with in the next section since the proof
in this case relies on the injections for the case $m\geq 1$.

Throughout this section, we assume that $m\geq 1$.
To construct an injection $\Phi$ from the set $Q(m,n)$ to the set $P(-m,n)$, we shall divide $Q(m,n)$ into six disjoint subsets $Q_i(m,n)$  $(1\leq i\leq 6)$
 and divide $P(-m,n)$ into eight disjoint subsets $P_{i}(-m,n)$ $(1\leq i\leq 8)$.  We then construct an injection $\Phi$ for $m\geq 1$
 consisting of six injections  $\phi_i$  from the set $Q_i(m,n)$  to the set $P_{i}(-m,n)$, where $1\leq i\leq 6$.

    Let $(\alpha,\beta)_{(m+j)\times j }$ be the $m$-Durfee rectangle symbol of a partition $\lambda$ in $Q(m,n)$. We write
    \[ \lambda=\left(\begin{array}{l}
    \alpha\\[3pt]
    \beta
    \end{array}\right)_{(m+j)\times j},\]
  and we say that the above $m$-Durfee rectangle symbol is a partition in $Q(m,n)$.
  In other words, $Q(m,n)$ is also considered as the set of
  $m$-Durfee rectangle symbols of partitions counted by $q(m,n)$.
    By Proposition \ref{pro-dufrec},
  we have either $j=0$ or $\beta_1=j$ with $j\geq 1$.   The subsets $Q_i(m,n)$ can  be described by using the $m$-Durfee
  rectangle symbol   $(\alpha,\beta)_{(m+j)\times j }$.
\begin{itemize}
\item[(1)] $Q_1(m,n)$ denotes the set of $m$-Durfee rectangle symbols
  in $Q(m,n)$
  for which either of the following conditions holds:\\
   (i) $j=0$;\\
   (ii) $j\geq 1$ and $\ell(\beta)-\ell(\alpha)\leq -1$; \\
   (iii)  $j\geq 1$, $\ell(\beta)-\ell(\alpha)=0$ and $\alpha_1=m+j$;

\item[(2)] $Q_2(m,n)$ denotes the set of $m$-Durfee rectangle symbols in $Q(m,n)$ such that $j\geq 1$, $\ell(\beta)-\ell(\alpha)\geq 0$ and $\alpha_1< m+j$;

\item[(3)] $Q_3(m,n)$ denotes the set of $m$-Durfee rectangle symbols in $Q(m,n)$   such that   $j\geq 1$,   $\ell(\beta)-\ell(\alpha)\geq 1$, $\alpha_1=m+j$ and $s(\beta)=1$;

\item[(4)] $Q_4(m,n)$ denotes the set of $m$-Durfee rectangle symbols in $Q(m,n)$ such that $j\geq 1$, $\ell(\beta)-\ell(\alpha)\geq 1$, $\alpha_1=m+j>\alpha_2$ and $s(\beta)\geq 2$;

\item[(5)] $Q_5(m,n)$ denotes the set of $m$-Durfee rectangle symbols in $Q(m,n)$ such that  $j\geq 1$, $\ell(\beta)-\ell(\alpha)\geq 1$, $\alpha_1=\alpha_2=m+j>\alpha_3$ and $s(\beta)\geq 2$;

\item[(6)] $Q_6(m,n)$ denotes the set of $m$-Durfee rectangle symbols in $Q(m,n)$ such that  $j\geq 1$, $\ell(\beta)-\ell(\alpha)\geq 1$, $\alpha_1=\alpha_2=\alpha_3=m+j$ and $s(\beta)\geq 2$.

\end{itemize}

To divide the set $P(-m,n)$ into eight disjoint subsets,
we also view $P(-m,n)$ as a set of $m$-Durfee rectangle symbols of partitions counted by $p(-m,n)$.   Let  $(\gamma,\delta)_{(m+j')\times j'}$ be the  $m$-Durfee rectangle symbol  in $P(-m,n)$.  By Proposition \ref{pro-dufrec}, we see that either  $j'=0$ or $j'\geq 1$ and $\ell(\delta)-\ell(\gamma)\leq 0$.   The subsets $P_i(-m,n)$  can  be described in terms of the $m$-Durfee   rectangle symbols   $(\gamma,\delta)_{(m+j')\times j'}$ in $P(-m,n)$:
\begin{itemize}
\item[(1)] $P_1(-m,n)$ denotes the set of  $m$-Durfee rectangle symbols   in $P(-m,n)$ for which either of the following conditions holds:\\
     (i) $j'=0$;\\
     (ii) $j'\geq 1$, $\ell(\delta)-\ell(\gamma)\leq -1$ and $\delta_1=j'$;\\
     (iii) $j'\geq 1$, $\ell(\gamma)=\ell(\delta)$, $\gamma_1=m+j'$ and $\delta_1=j'$;

\item[(2)] $P_2(-m,n)$ denotes the set of  $m$-Durfee rectangle symbols  in $P(-m,n)$  such that $j'\geq 1$  and $\delta_1=j'-1$;

\item[(3)] $P_3(-m,n)$ denotes the set of  $m$-Durfee rectangle symbols in $P(-m,n)$ such that $j'\geq 2$  and $\delta_1\leq j'-2$;

\item[(4)] $P_4(-m,n)$ denotes the set of $m$-Durfee rectangle symbols in $P(-m,n)$ such that $j'\geq 1$, $\ell(\gamma)=\ell(\delta)$, $\gamma_1=m+j'-1$, $\delta_1=j'$ and $\delta$ has a part equal to $2$;

\item[(5)] $P_5(-m,n)$ denotes the set of $m$-Durfee rectangle symbols in $P(-m,n)$ such that $j'\geq 1$, $\ell(\gamma)=\ell(\delta)$, $\gamma_1\leq m+j'-3$ and $\delta_1=j'$;

\item[(6)] $P_6(-m,n)$ denotes the set of $m$-Durfee rectangle symbols in $P(-m,n)$ such that $j'\geq 1$, $\ell(\gamma)=\ell(\delta)$, $\gamma_1= m+j'-2$ and $\delta_1=j'$;

\item[(7)] $P_7(-m,n)$ denotes the set of $m$-Durfee rectangle symbols in $P(-m,n)$ such that $j'\geq 1$, $\ell(\gamma)=\ell(\delta)$, $\gamma_1=m+j'-1>\gamma_2$, $\delta_1=j'$ and $\delta$ has no parts equal to $2$;

\item[(8)] $P_8(-m,n)$ denotes the set of $m$-Durfee rectangle symbols in $P(-m,n)$ such that $j'\geq 1$, $\ell(\gamma)=\ell(\delta)$, $\gamma_1=\gamma_2=m+j'-1$, $\delta_1=j'$ and $\delta$ has no parts equal to $2$.

\end{itemize}

We are now ready to present six injections $\phi_i$ from the set $Q_i(m,n)$ to the set $P_{i}(-m,n)$, where $1\leq i\leq 6$.
It is clear  that  $Q_1(m,n)$ coincides with $P_1(-m,n)$. So
$\phi_1$ can be set to be the identity map. The following
lemma gives an injection from $Q_2(m,n)$ to $P_2(-m, n)$.

 \begin{lem}\label{phi-2} For $m\geq 0$,
there is an injection $\phi_2$ from the set $Q_2(m,n)$ to the set $P_2(-m,n)$.
\end{lem}

\pf  Let
 \[\lambda=\left(\begin{array}{c }
\alpha\\[3pt]
\beta
\end{array}
\right)_{(m+j)\times j}=\left(\begin{array}{cccc }
\alpha_1,&\alpha_2,&\ldots, & \alpha_s\\[3pt]
\beta_1,& \beta_2,&\ldots, & \beta_t
\end{array}
\right)_{(m+j)\times j} \]
be an $m$-Durfee rectangle symbol in $Q_2(m,n)$.
By definition, we have    $\beta_1=j\geq 1$, $\alpha_1<m+j$ and $t-s\geq 0$.

Define
   \[\phi_2(\lambda)=\left(\begin{array}{c }
\gamma\\[3pt]
\delta
\end{array}
\right)_{(m+j')\times j'}=\left(\begin{array}{cccccc }
\alpha_1+1,&\alpha_2+1,&\ldots, &\alpha_s+1,& 1^{t-s}\\[3pt]
\beta_1-1,& \beta_2-1,&\ldots, &\beta_t-1&
\end{array}
\right)_{(m+j)\times j}.\]
It is evident  that    $\ell(\delta)\leq t$ and $\ell(\gamma)=t$. So $\ell(\delta)-\ell(\gamma)\leq 0$. Moreover it is easy to see that $\delta_1=j-1$ and  $|\phi_2(\lambda)|=|\lambda|$. This proves that $\phi_2(\lambda)$ is in $ P_2(-m,n).$

To prove that the map $\phi_2$ is an injection, let
\[H(m,n)=\{\phi_2(\lambda)\colon \lambda \in Q_2(m,n)\}.\]
If $n\neq m+1$, it is easy to check that $H(m,n)$ is coincide with $P_2(-m,n)$; If $n=m+1$, we have
\[H(m,n)=P_2(-m,n)\setminus \{(\emptyset,\emptyset)_{(m+1)\times 1}\}.\]

Let
\[\mu= \left(\begin{array}{cc}
 \gamma\\[3pt]
\delta
\end{array}\right)_{(m+j')\times j'}=\left(\begin{array}{ccccccc}
\gamma_1,&\gamma_2,&\ldots, &  \gamma_{s'}\\[3pt]
\delta_1,& \delta_2,&\ldots, & \delta_{t'}
\end{array}
\right)_{(m+j')\times j'}\]
 be an $m$-Durfee rectangle symbol in $H(m,n)$. Since $\mu \in P_2(-m,n)$, we have $s'\geq t'$. Define $\sigma(\mu)$ to be
 \[\sigma(\mu)=\left(\begin{array}{ccccccc}
\gamma_1-1,&\gamma_2-1,&\ldots, &  \gamma_{s'}-1\\[3pt]
\delta_1+1,& \delta_2+1,&\ldots, & \delta_{t'}+1, & 1^{s'-t'}
\end{array}
\right)_{(m+j')\times j'}.
 \]
It can be checked $\sigma(\mu)$ is in $Q_2(m,n)$ and  $\sigma(\phi_2(\lambda))=\lambda$ for any $\lambda$ in $Q_2(m,n)$. Hence the map $\phi_2$  is a bijection between $Q_2(m,n)$ and $H(m,n)$. \qed

For example, for $m=2$ and $n=31$, let
  \[\lambda=\left(\begin{array}{cccc }
4,&2,&2\\[3pt]
3,&2,&2,&1
\end{array}
\right)_{5\times 3}\]
be a $2$-Durfee rectangle symbol in $ Q_2(2,31)$.
Applying the map $\phi_2$ to $\lambda$, we obtain
 \[\phi_2(\lambda)=\left(\begin{array}{cccc }
5,&3,&3,&1\\[3pt]
2,&1,&1&
\end{array}
\right)_{5\times 3},\]
which is a $2$-Durfee rectangle symbol belonging to $P_2(-2,31)$. Applying $\sigma$ to $\phi_2(\lambda)$, we recover $\lambda$, that is, $\sigma(\phi_2(\lambda))=\lambda$.

 \begin{lem}\label{phi-3}For $m\geq 0$, there is a bijection $\phi_3$ between the set $Q_3(m,n)$ and the set $P_3(-m,n)$.
\end{lem}

 \pf Let
 \[\lambda=\left(\begin{array}{c }
\alpha\\[3pt]
\beta
\end{array}
\right)_{(m+j)\times j}=\left(\begin{array}{cccc }
\alpha_1,&\alpha_2,&\ldots, & \alpha_s\\[3pt]
\beta_1,& \beta_2,&\ldots, & \beta_t
\end{array}
\right)_{(m+j)\times j}\]
be an $m$-Durfee rectangle symbol in $Q_3(m,n)$. By definition, we see that  $j=\beta_1\geq \beta_t=1$, $\alpha_1=m+j$ and $t-s\geq 1$.

Define
    \[\phi_3(\lambda)=\left(\begin{array}{c }
\gamma\\[3pt]
\delta
\end{array}
\right)_{(m+j')\times j'}=\left(\begin{array}{cccc }
\alpha_2+1,&\ldots ,& \alpha_s+1, &1^{t-s-1}\\[3pt]
  \beta_2-1,&\ldots ,& \beta_{t}-1&
\end{array}
\right)_{(m+j+1)\times (j+1)}.\]
 To prove that $\phi_3(\lambda) \in P_3(-m,n)$, we need to verify that $\gamma_1\leq m+j'$, $\delta_1\leq j'-2$, $\ell(\delta)-\ell(\gamma) \leq 0$  and $|\lambda|=|\phi_3(\lambda)|$.
First, it is easy to see that
\[\gamma_1=\alpha_2+1\leq m+j+1=m+j'\]
and
\[\delta_1=\beta_2-1\leq j-1\leq j'-2.\]
 By definition, we have $\ell(\gamma)=t-2$ and $\ell(\delta)\leq t-2$ for $\beta_t=1$. Hence $\ell(\delta)-\ell(\gamma) \leq 0$.

Note that
\begin{align*}
|\phi_3(\lambda)|&=|\gamma|+|\delta|+(j+1)(m+j+1) .
\end{align*}
But
 \begin{eqnarray*}
|\gamma|+|\delta|&=&(|\alpha|-\alpha_1+t-2)+(|\beta|-\beta_1-(t-1))\\[3pt]
&=&|\alpha|+|\beta|-(m+j)-j-1,
\end{eqnarray*}
we find that
\begin{align*}
|\phi_3(\lambda)|&=|\alpha|+|\beta|-(m+j)-1-j+(j+1)(m+j+1) \\[3pt]
&=|\alpha|+|\beta|+j(m+j)  ,
\end{align*}
which equals $|\lambda|$.
Hence $\phi_3(\lambda)\in P_3(-m,n)$. Moreover, it can be checked that
   $\phi_3$ is invertible. So we conclude that $\phi_3$ is a  bijection. \qed

For example, for $m=2$ and $n=34$,  let
\[\lambda= \left(\begin{array}{cccc }
5,&4,&1\\[3pt]
3,&3,&2,&1
\end{array}
\right)_{5\times 3}
\]
be a $2$-Durfee rectangle symbol in $Q_3(2,34)$.
Applying the bijection $\phi_3$ to $\lambda$,  we get
\[\phi_3(\lambda)=\left(\begin{array}{cccc }
5,&2\\[3pt]
2,&1
\end{array}
\right)_{6\times 4},\]
which is in  $ P_3(-2,34)$.

The following proposition will be used in the construction of the injection $\phi_4$.

\begin{pro}\label{pro-phi-4}
For $m\geq 0$, let
\[\lambda=\left(\begin{array}{c }
\alpha\\[3pt]
\beta
\end{array}
\right)_{(m+j)\times j}
= \left(\begin{array}{cccc }
\alpha_1,&\alpha_2,&\ldots, & \alpha_s\\[3pt]
\beta_1,& \beta_2,&\ldots, & \beta_t
\end{array}
\right)_{(m+j)\times j} \]
be an $m$-Durfee rectangle symbol in $Q_4(m,n)$.
Then there exists an integer $1\leq k\leq s$ such that
\begin{equation}\label{pro-r-1}
{\alpha}_{k+1}\leq  {\beta}_k-1
\end{equation}
and
\begin{equation}\label{pro-r-2}
{\alpha}_{k}\geq {\beta}_{k+1}-1.
\end{equation}
 \end{pro}

\pf   By definition of $Q_4(m,n)$, we have   $j=\beta_1\geq \beta_t\geq 2$, $m+j=\alpha_1>\alpha_2$ and $t-s\geq 1$. When $m=0$, we may choose $k=1$, since
\[ {\alpha}_2\leq j-1= {\beta}_1-1\]
and
\[\alpha_{1}=j>  {\beta}_{2}-1.\]
 When $m\geq 1$, let
\[h=\min\{i \colon 1\leq i\leq t, {\alpha}_i\leq  {\beta}_i-1\}.\]
Setting $k=h-1$, we  proceed to show that $1\leq k\leq s$ and relations (\ref{pro-r-1}) and (\ref{pro-r-2}) hold. Since $ {\beta}_t\geq 2$,
$ {\alpha}_{s+1}=0$ and $t\geq s+1$,  we have $ {\alpha}_{s+1}\leq {\beta}_{s+1}-1$, which implies that $h\leq s+1$, that is, $k\leq s$. Observing that $\alpha_1=j+m> j-1=\beta_1-1$, we get $h\geq 2$, that is, $k\geq 1$. Thus, we have $1\leq k\leq s$. By the definition of $h$,  we find that
 \begin{equation*}
  {\alpha}_h\leq {\beta}_h-1
 \end{equation*}
 and
 \begin{equation*}
  {\alpha}_{h-1}>  {\beta}_{h-1}-1.
\end{equation*}
It follows that
 \[  {\alpha}_{h}\leq  {\beta}_{h}-1\leq {\beta}_{h-1}-1 \]
 and
\[ {\alpha}_{h-1}>  {\beta}_{h-1}-1\geq  {\beta}_{h}-1,\]
which implies that $k=h-1$. This completes the proof. \qed

 \begin{lem}\label{phi-4}For $m\geq 0$, there is an injection $\phi_4$ from the set $Q_4(m,n)$ to the set $P_4(-m,n)$.
\end{lem}

\pf We construct a map $\phi_4$ from the set $Q_4(m,n)$ to the set $P_4(-m,n)$,
      then we show that it is an injection. Let
\[\lambda=\left(\begin{array}{c }
\alpha\\[3pt]
\beta
\end{array}
\right)_{(m+j)\times j}
= \left(\begin{array}{cccc }
\alpha_1,&\alpha_2,&\ldots, & \alpha_s\\[3pt]
\beta_1,& \beta_2,&\ldots, & \beta_t
\end{array}
\right)_{(m+j)\times j} \]
be an $m$-Durfee rectangle symbol in $Q_4(m,n)$.
By Proposition \ref{pro-phi-4}, we may choose $k$ to be  the minimum integer
such that $1\leq k \leq s$, $\alpha_{k+1}\leq  {\beta}_k-1$ and ${\alpha}_{k}\geq {\beta}_{k+1}-1$. By the definition of $Q_4(m,n)$, we have   $j=\beta_1\geq \beta_t\geq 2$, $m+j=\alpha_1>\alpha_2$ and $t-s\geq 1$. So we may define
\begin{eqnarray}  \phi_4(\lambda)&=& \left(\begin{array}{cc}
 \gamma\\[3pt]
\delta
\end{array}\right)_{(m+j')\times j'} \nonumber \\[5pt] \label{phi-4-temp}
&=&\left(\begin{array}{cccccccccc}
 {\alpha}_1-1,& {\alpha}_2,&\ldots, &  {\alpha}_k,& {\beta}_{k+1}-1, &\ldots, & {\beta}_t-1\\[3pt]
 {\beta}_1,&  {\beta}_2,&\ldots, &  {\beta}_k, & {\alpha}_{k+1}+1, &\ldots, & \alpha_s+1, &2,& 1^{t-s-1}
\end{array}
\right)_{(m+j)\times j}.
\end{eqnarray}
 Apparently,  $\gamma_1= {\alpha}_1-1=j'+m-1$, $\delta_1= {\beta}_1=j'$, $\ell(\gamma)=\ell(\delta)=t$ and  $\delta_{s+1}=2$. Furthermore, it can be easily checked that $|\phi_4(\lambda)|=|\lambda|$.
This yields that $\phi_4(\lambda)  \in P_4(-m,n)$.

To prove that $\phi_4$ is an injection, Let
\begin{equation*}\label{phi4-im}
I(m,n)=\{\phi_4(\lambda)\colon \lambda\in Q_4(m,n)\}
\end{equation*}
be the set of images of $\phi_4$,
which has been shown to be a subset of $P_4(-m,n)$.
We wish to show that the construction of $\phi_4$ is
reversible, which implies that $\phi_4$ is an injection.
More precisely, we shall show that there exists a map
$\varphi$ from $I(m,n)$ to $Q_4(m,n)$ such that
for any $\lambda$ in $Q_4(m,n)$ we have
\begin{equation*}\label{r-phi-varphi}
\varphi(\phi_4(\lambda))=\lambda.
\end{equation*}

We now describe the map $\varphi$.
Let
\begin{equation}\label{defi-mu} \mu= \left(\begin{array}{cc}
 \gamma\\[3pt]
\delta
\end{array}\right)_{(m+j')\times j'}=\left(\begin{array}{ccccccc}
\gamma_1,&\gamma_2,&\ldots, &  \gamma_{t'}\\[3pt]
\delta_1,& \delta_2,&\ldots, & \delta_{t'}
\end{array}
\right)_{(m+j')\times j'}
\end{equation}
be an $m$-Durfee rectangle symbol in $I(m,n)$. The following procedure
generates an $m$-Durfee rectangle symbol $\varphi(\mu)$ that is in $Q_4(m,n)$.

We claim that for $\mu \in I(m,n)$ given by \eqref{defi-mu}, there exists an
 integer $k'$ such that
$1\leq k' \leq \ell(\gamma)-1$ and
\begin{equation}\label{cond-varphi}
\delta_{k'}-1\ge \gamma_{k'+1}, \quad \gamma_{k'}\ge \delta_{k'+1}-1\geq 1.
\end{equation}
Since $\mu \in I(m,n)$,   there exists $\lambda \in Q_4(m,n)$ such that $\phi_4(\lambda)=\mu$.    By the choice of $k$
in the construction $\phi_4(\lambda)$, we see that
   \begin{equation*}\label{phi-4-k-stri}
   1\leq k\leq s\leq t-1=\ell(\gamma)-1.
   \end{equation*}
   Again, from the construction
   \eqref{phi-4-temp} of $\phi_4(\lambda)$, we find that
    \begin{equation*}\label{phi-4-a}
    \delta_k \ge \gamma_{k+1}+1
    \end{equation*}
and
\begin{equation*}\label{phi-4-b}
\gamma_k \ge \delta_{k+1}-1\ge 1.
\end{equation*}
So $k$ satisfies the conditions in \eqref{cond-varphi}.
Thus the claim is verified.

Now, we may choose  $k'$ to be the minimum integer such that  $1\leq k' \leq \ell(\gamma)-1,\, \delta_{k'}-1\ge \gamma_{k'+1}$ and $\gamma_{k'}\ge \delta_{k'+1}-1\geq 1$.  Since
 $\mu$ is in $P_4(-m,n)$,  the partition  $\delta$
  in the $m$-Durfee rectangle symbol of $\mu$
  has a part equal to $2$. Assume that $\delta_{s'}=2>\delta_{s'+1}$.
  Hence we may define
\begin{eqnarray}
\varphi(\mu)&=&\left(\begin{array}{cc}
 { \alpha}\\[3pt]
 {\beta}
\end{array}\right)_{(m+j)\times j} \nonumber \\[5pt] \label{constr}
&=&\left(\begin{array}{ccccccc}
\gamma_1+1,&\gamma_2,&\ldots,   &\gamma_{k'},&\delta_{k'+1}-1, &\ldots, &\delta_{s'-1}-1\\[3pt]
\delta_1,&\delta_2,&\ldots, &\delta_{k'},&\gamma_{k'+1}+1, &\ldots, &\gamma_{t'}+1
\end{array}
\right)_{(m+j')\times j'}.
\end{eqnarray}
Evidently, $ {\beta}_1=\delta_1=j$, $ {\alpha}_1=\gamma_1+1=m+j> {\alpha}_2$, $ {\beta}_{t'}=\gamma_{t'}+1\geq 2$ and $t'>s'-1$. Moreover, it is easy to check that $|\varphi(\mu)|=|\mu|$. So we deduce that $\varphi(\mu) \in Q_4(m,n).$

It remains to verify that  $\varphi(\phi_4(\lambda))=\lambda$.
By the constructions \eqref{phi-4-temp} and \eqref{constr} of $\phi_4(\lambda)$ and $\varphi(\mu)$,   it suffices  to show that the integer $k$ appearing
 in the representation of $\phi_4(\lambda)$ coincides with the
 integer $k'$ appearing in
the representation of $\varphi(\phi_4(\lambda))$.

Recall that  $k$ is the minimum integer determined by  $\lambda$  subject to the conditions
\begin{equation}\label{cond-k}
1\leq k \leq s, \quad  {\alpha}_{k}\geq {\beta}_{k+1}-1, \quad  \text{and} \quad \alpha_{k+1}\leq  {\beta}_k-1.
\end{equation}
On the other hand, it can be shown that $k$ is also the minimum integer $k'$ depending on
 $\phi_4(\lambda)$   such that
 \begin{equation}\label{cond-k'}
1\leq k' \leq \ell(\gamma)-1, \quad \delta_{k'}-1\ge \gamma_{k'+1} \quad  \text{and} \quad \gamma_{k'}\ge \delta_{k'+1}-1\geq 1.
\end{equation}
 From the definitions of $k$ and $s$, we find that $s\leq t-1=\ell(\gamma)-1$, which implies $k \leq \ell(\gamma)-1$.
 By the construction \eqref{phi-4-temp} in $\phi_4(\lambda)$,  we have
$\gamma_{k+1}=\beta_{k+1}-1$ and $\delta_k=\beta_k.$
Furthermore, we have
$\gamma_1=\alpha_1-1$, and $\gamma_k=\alpha_k$  when $k\geq 2$. It can also be seen that $\delta_{s+1}=2$ and $\delta_{k+1}=\alpha_{k+1}+1$ when $1\leq k\leq s-1$.
Hence we deduce that  $\delta_{k}-1\ge \gamma_{k+1}$ and
 $\gamma_{k}\ge \delta_{k+1}-1\geq 1$ for $1\leq k \leq s$.
In other words, $k$ satisfies the conditions in \eqref{cond-k'}.

Finally, we need to show that $k$ is   the minimum integer satisfying conditions in \eqref{cond-k'}.  Assume to the contrary that there is an
 integer  $1\leq p\leq k-1$  for which the conditions in \eqref{cond-k'}
 are satisfied, that is,
\[\delta_{p}-1\ge \gamma_{p+1} \quad  \text{and} \quad \gamma_{p}\ge \delta_{p+1}-1\geq 1.\]
 From construction \eqref{phi-4-temp} of $\phi_4(\lambda)$ and the
 assumption $1\leq p \leq k-1$, we find that
 \[\alpha_{p+1}=\gamma_{p+1},\quad \beta_p=\delta_p,\quad \beta_{p+1}=\delta_{p+1}.\]
Moreover, by (\ref{phi-4-temp}) we see that   $\alpha_p=\gamma_p+1$ if
 $p=1$ and $\alpha_p=\gamma_p$ if $p\geq 2$.
In either case, we have
 \[{\alpha}_{p}\geq {\beta}_{p+1}-1 \quad  \text{and} \quad \alpha_{p+1}\leq  {\beta}_p-1.\]
 This means that $p$ also satisfies the conditions in \eqref{cond-k},
 contradicting the minimality of $k$. So we conclude that
 $k$ is  the minimum integer satisfying conditions in \eqref{cond-k'},
 which implies that $\varphi(\phi_4(\lambda))=\lambda$.
   This completes the proof. \qed

For example, for $m=2$ and $n=41$, consider the following
 $2$-Durfee rectangle symbol in $Q_4(2,41)$:
\[\lambda=\left(\begin{array}{ccccccc }
5,&4,&2,&1\\[3pt]
3,&3,&2,&2,&2,&2
\end{array}
\right)_{5\times 3} .\]
It can be checked that $k=2$. Applying the injection $\phi_4$ to $\lambda$, we have
\[\mu=\phi_4(\lambda)=\left(\begin{array}{ccccccc }
4,&4,&1,&1,&1,&1\\[3pt]
3,&3,&3,&2,&2,&1
\end{array}
\right)_{5\times 3},\]
which is in $P_4(-2,41)$.
Applying $\varphi$ to $\mu$, we have  $k'=2$ and
 $\varphi(\mu)=\lambda$.

We next describe the injection $\phi_5$ from $Q_5(m,n)$ to $P_5(-m,n)$ for the case $m\geq 1$. It should be noted that
 the construction for $m\geq 1$ does
 not apply to the case $m=0$.

\begin{lem}\label{phi-5}For $m\geq 1$, there is an injection $\phi_5$ from the set $Q_5(m,n)$ to the set $P_5(-m,n)$.
\end{lem}

\pf Let
\begin{equation*}\label{lambda-phi5}
\lambda=\left(\begin{array}{c }
\alpha\\[3pt]
\beta
\end{array}
\right)_{(m+j)\times j}
= \left(\begin{array}{cccc }
\alpha_1,&\alpha_2,&\ldots, & \alpha_s\\[3pt]
\beta_1,& \beta_2,&\ldots, & \beta_t
\end{array}
\right)_{(m+j)\times j}
\end{equation*}
be an $m$-Durfee rectangle symbol in $Q_5(m,n)$.  By definition, we have $j=\beta_1\geq \beta_t \geq 2$, $\alpha_1=\alpha_2=m+j>\alpha_3$ and $t-s\geq 1$.

 Since $\alpha_2-m+2=j+2 > \beta_3-1$, we may  choose the maximum number $k$ such that $1\leq k\leq t-1$ and $\alpha_k-m+2\ge \beta_{k+1}-1$.
 To define $\phi_5(\lambda)$, we construct two partitions $\gamma$ and $\delta$.  It is clear that $k\geq 2$.
So we may define
\begin{equation}\label{gamma-1}
\gamma=(\beta_2+m-2,\,\ldots,\, \beta_k+m-2,\, \alpha_{k+1}+1,\,\ldots,\,\alpha_{t}+1)
\end{equation}
and
\begin{equation}\label{delta-1a}
\delta=(\alpha_2+1-m,\,\alpha_3+2-m,\,\ldots, \, \alpha_k+2-m,\,\beta_{k+1}-1,\,\ldots,\,\beta_{t}-1).
\end{equation}
Notice that when $k=2$ the above definition \eqref{delta-1a}
may be endangered by ambiguity.
In this case, \eqref{delta-1a} is interpreted as
\begin{equation*}\label{delta-1b}
\delta=(\alpha_2+1-m,\,\beta_{3}-1,\,\ldots,\,\beta_{t}-1).
\end{equation*}
We now define
 \begin{equation}\label{phi-5-lambda}\phi_5(\lambda)=\left(\begin{array}{cc}
 \gamma\\[3pt]
\delta
\end{array}\right)_{(m+j+1)\times (j+1)}.
\end{equation}

We first prove that $(\gamma,\delta)_{(m+j+1)\times (j+1)}$
is an $m$-Durfee rectangle symbol. To this end, we need to show that $\gamma$ and $\delta$ are partitions with $\gamma_1\leq m+j+1$ and $\delta_1\leq j+1$.  We then verify that $(\gamma,\delta)_{(m+j+1)\times (j+1)}$ satisfies the
conditions for $P_5(-m,n)$.

 To prove that $\delta$ is a partition, it suffices to show that when $k=2$, we have
\begin{equation}\label{delta-k-2}
\alpha_2+1-m\geq \beta_{3}-1,
\end{equation}
and when $k\geq 3$, we have
\begin{equation}\label{delta-k-b}
\alpha_2+1-m\geq \alpha_3+2-m
\end{equation}
and
\begin{equation}\label{delta-k-a}
\alpha_k+2-m\geq \beta_{k+1}-1.
\end{equation}

 When $k=2$, since $\alpha_2-m+1=j+1$ and $\beta_3\leq \beta_1=j$,  we see that \eqref{delta-k-2} holds, and so  $\delta$ is a partition. When $k\geq 3$, since $\alpha_2>\alpha_3$, we get \eqref{delta-k-b}.  On the other hand, \eqref{delta-k-a}  follows from the choice of $k$.  Hence  $\delta$ forms  a partition when $k\geq 3$.

We now verify that $\gamma$ is a partition. From the definition \eqref{gamma-1} of $\gamma$, it suffices to show that
\begin{equation}\label{gamma}
\beta_k+m-2\geq  \alpha_{k+1}+1.
\end{equation}
Keep in mind that $k$ is in the range from $2$ to $t-1$.
When $k=t-1$,   \eqref{gamma} becomes $\beta_{t-1}+m-2\geq \alpha_t+1$, which is valid since $\beta_{t-1}\geq 2$ and $\alpha_t=0$. When $2\leq k\leq t-2$,  since $k$ is the maximum integer such that
$\alpha_k-m+2\geq \beta_{k+1}-1,$ we find $\alpha_{k+1}-m+2< \beta_{k+2}-1$, which implies \eqref{gamma}.
This proves that $\gamma$ is a partition.

Next we demonstrate that $(\gamma,\delta)_{(m+j+1)\times (j+1)}$
is an $m$-Durfee rectangle symbol in $P_5(-m,n)$.
It is clear from \eqref{gamma-1} and \eqref{delta-1a}  that
  $\delta_1=\alpha_2+1-m=j+1$, $\gamma_1=\beta_2+m-2\leq j+m-2$, and  $\ell(\gamma)=\ell(\delta)=t-1$. It remains to check that $|(\gamma,\delta)_{(m+j+1)\times (j+1)}|=|\lambda|$.
Note that
\begin{eqnarray*}
|\gamma|+|\delta|&=&
| {\alpha}|- {\alpha}_1+(2-m)  (k-2)+1-m+(t-k)\\[3pt]
&&\quad +| {\beta}|- {\beta}_1+(m-2) (k-1)-(t-k)\\[3pt]
&=&| {\alpha}|- {\alpha}_1+| {\beta}|- {\beta}_1-1.
\end{eqnarray*}
Using $\beta_1=j$ and $\alpha_1=m+j$, so we get
\[|\gamma|+|\delta|=|\alpha|+|\beta|-(2j+m+1).
\]
Hence
\begin{eqnarray*}
|(\gamma,\delta)_{(m+j+1)  \times (j+1)}|&=&|\gamma|+|\delta|+(m+j+1)  (j+1)\\[3pt]
&=&|\alpha|+|\beta|+j (j+m),
\end{eqnarray*}
which equals $|\lambda|$. So we arrive at the
conclusion $(\gamma,\delta)_{(m+j+1)\times (j+1)}\in P_5(-m,n)$.

Next we proceed to prove that  $\phi_5$ is an injection. Let
\begin{equation*}
J(m,n)=\{\phi_5(\lambda)\colon \lambda\in Q_5(m,n)\}
\end{equation*}
be the set of images of $\phi_5$. It has been shown
that $J(m,n)$ is a subset of $P_5(-m,n)$.
We wish to construct a map
$\tau$ from $J(m,n)$ to $Q_5(m,n)$ such that
for any $\lambda$ in $Q_5(m,n)$, we have
\begin{equation*}\label{r-phi5-varphi}
\tau(\phi_5(\lambda))=\lambda.
\end{equation*}

To describe the map $\tau$, let
\begin{equation*}\label{defi-phi-5-mu} \mu= \left(\begin{array}{cc}
 \gamma\\[3pt]
\delta
\end{array}\right)_{(m+j')\times j'}=\left(\begin{array}{ccccccc}
\gamma_1,&\gamma_2,&\ldots, &  \gamma_{t'}\\[3pt]
\delta_1,& \delta_2,&\ldots, & \delta_{t'}
\end{array}
\right)_{(m+j')\times j'}
\end{equation*}
be an $m$-Durfee rectangle symbol in $J(m,n)$, that is,
    there is an $m$-Durfee rectangle symbol $\lambda=(\alpha,\beta)_{(m+j)\times j}$ in $Q_5(m,n)$ such that
$\phi_5(\lambda)=\mu$.
We claim that $\gamma_{t'}=1$  and there exists an integer $k'$ such that
\begin{equation}\label{cond-5-k'}
1\leq k'\leq t'-1\quad \text{and} \quad \gamma_{k'}-m+1\ge\delta_{k'+1}.
\end{equation}
From the constructions  \eqref{gamma-1} and \eqref{delta-1a} of $\phi_5$, we see that $\gamma_{t'}=\alpha_{t}+1=1$,  $\gamma_{k-1}=\beta_k+m-2$ and $\delta_k=\beta_{k+1}-1$.
It follows that $\gamma_{k-1}-m+1\ge\delta_{k}$.  Since $1\leq k-1\leq t-2=t'-1$, we reach the conclusion that $k-1$ satisfies the conditions in \eqref{cond-5-k'}. This proves the claim.

Now we may choose $k'$ to be the maximum integer such that $1\leq k'\leq t'-1$ and
\begin{equation}\label{phi-5-i-t-1}
\gamma_{k'}-m+1\geq \delta_{k'+1}.
\end{equation}
The choice of $k'$ yields that $\gamma_{k'+1}-m+1<\delta_{k'+2}$ when $1\leq k'\leq t'-2$, which implies $\gamma_{k'+1}-1<\delta_{k'}-2+m$.  When $k'=t'-1$,
we also have $\gamma_{k'+1}-1\leq \delta_{k'}-2+m$
since  $\gamma_{t'}=1$. Combining the above two cases for $k'$,  we find that
\begin{equation}\label{phi-5-i-t-2}
\gamma_{k'+1}-1\leq \delta_{k'}-2+m.
\end{equation}

By \eqref{phi-5-i-t-1} and \eqref{phi-5-i-t-2}, we may define
\begin{equation*}
\tau(\mu)=\left(\begin{array}{cc}
{\alpha}\\[3pt]
{\beta}
\end{array}\right)_{(m+j'-1)\times (j'-1)},
\end{equation*}
where
\begin{equation}
\alpha=(j'+m-1,\,\delta_1-1+m,\,\delta_2-2+m,\,\ldots, \,\delta_{k'}-2+m,\,\gamma_{k'+1}-1,\,\ldots,\, \gamma_{t'}-1)
\end{equation}
and
\begin{equation}
\beta=(j'-1,\,\gamma_1+2-m,\,\ldots,\,\gamma_{k'}+2-m,\,\delta_{k'+1}+1,\, \ldots,\,\delta_{t'}+1).
\end{equation}
 It is easily checked that $\tau(\mu)\in Q_5(m,n)$.

We are now ready to verify that  $\tau(\phi_5(\lambda))=\lambda$.
By the constructions of $\phi_5(\lambda)$ and $\tau(\mu)$,   it suffices  to show that the integer $k$ appearing
 in the representation of $\phi_5(\lambda)$ is equal to the
 integer $k'$ appearing in
the representation of $\tau(\phi_5(\lambda))$ plus $1$, namely, $k'=k-1.$

Recall that  $k$ is the maximum integer determined by  $\lambda$  subject to the conditions
\begin{equation}\label{cond-k-5}
1\leq k\leq t-1\quad\text{and}\quad\alpha_k-m+2\ge \beta_{k+1}-1.
\end{equation}
On the other hand, it can be shown that $k-1$ is   the maximum integer $k'$ determined by $\phi_5(\lambda)$  subject to the
conditions
 \begin{equation}\label{cond-k'-5}
1\leq k'\leq t'-1 \quad \text{and}\quad \gamma_{k'}-m+1\geq \delta_{k'+1}.
\end{equation}
It is not difficult to check that $k-1$ satisfies the conditions in \eqref{cond-k'-5}.
It remains to show that $k-1$ is   the maximum integer satisfying the conditions in \eqref{cond-k'-5}.  Assume to the contrary that there is an
 integer  $p\geq k$  for which the conditions in \eqref{cond-k'-5}
 are satisfied, that is, $k\leq p\leq t'-1$ and
 \begin{equation}\label{gamma-p}
  \gamma_{p}-m+1\geq \delta_{p+1}.
  \end{equation}
Since $t'=t-1$, we have
\begin{equation}\label{kpt2}
k\leq p\leq t-2.
\end{equation}From constructions \eqref{gamma-1} and \eqref{delta-1a} of $\phi_5(\lambda)$, we find that
$\gamma_{p}=\alpha_{p+1}+1$ and $\delta_{p+1}=\beta_{p+2}-1.$
By \eqref{gamma-p},  we deduce that ${\alpha}_{p+1}-m+2\geq {\beta}_{p+2}-1$. Moreover, by (\ref{kpt2}) we get $k+1\leq p+1\leq t-1$. Thus,  \eqref{cond-k-5}
is valid with $k$ being replaced by $p+1$. But this
contradicts the choice of $k$.  So we conclude that
 $k-1$ is  the maximum integer satisfying conditions in \eqref{cond-k'-5},
 which implies that $\tau(\phi_5(\lambda))=\lambda$. This completes the proof. \qed

For example, for $m=1$ and $n=34$,  consider the following
 $1$-Durfee rectangle symbol in $Q_5(1,34)$:
\[\lambda=\left(\begin{array}{ccccccc }
4,&4,&2\\[3pt]
3,&3,&2,&2,&2
\end{array}
\right)_{4\times 3}.\]
It can be checked that $k=4$.  Applying the injection $\phi_5$ to $\lambda$,  we have
\[\mu=\phi_5(\lambda)=\left(\begin{array}{ccccccc }
2,&1,&1,&1\\[3pt]
4,&3,&1,&1
\end{array}
\right)_{5\times 4},\]
which is in $P_5(-1,34)$. Applying $\tau$ to $\mu$, we have $k'=3$ and $\tau(\mu)=\lambda$.

It should be remarked that the injection $\phi_5$ is not valid for the case $m=0$. More precisely,  $\phi_5$ does not apply to Durfee  symbols $\lambda=(\alpha,\beta)_j$ in $Q_5(0,n)$ with $\beta_{t-1}=2$ where $\ell(\beta)=t$ and $\ell(\alpha)=s<t$. Assume that $\beta_{t-1}=2$, then we have $\alpha_{t-1}+2\ge 2>\beta_{t}-1$, and hence  $k=t-1$. Applying $\phi_5$ to $(\alpha,\beta)_j$, we
get
\[
\gamma=(\beta_2-2,\,\ldots,\, \beta_{t-1}-2,\, \alpha_{t}+1),
\]
which is not a partition, since $\gamma_{t-2}=\beta_{t-1}-2=0$ and $\gamma_{t-1}=\alpha_t+1=1$.

  In the following lemma, we give
an injection $\phi_6$ from $Q_6(m,n)$ to $P_6(-m,n)$ for the case $m\geq 1$.
Notice that this injection is not valid for $m=0$.

 \begin{lem}\label{phi-6}For $m\geq 1$, there is an injection $\phi_6$ from  the set $Q_6(m,n)$ to the set $P_6(-m,n)$.
\end{lem}

\pf To define the map $\phi_6$, let
\begin{equation*} \label{lambda-phi6}
\lambda=\left(\begin{array}{c }
\alpha\\[3pt]
\beta
\end{array}
\right)_{(m+j)\times j}= \left(\begin{array}{cccc }
\alpha_1,&\alpha_2,&\ldots, & \alpha_s\\[3pt]
\beta_1,& \beta_2,&\ldots, & \beta_t
\end{array}
\right)_{(m+j)\times j}
\end{equation*}
be an $m$-Durfee rectangle symbol in $Q_6(m,n)$. By definition, we have $j=\beta_1\geq \beta_t \geq 2$, $\alpha_1=\alpha_2=\alpha_3=m+j$ and $t-s\geq 1$.

Since $\alpha_3-m+1=j+1>\beta_3-1$,
there exists the maximum integer $k$ such that $  k\le s$ and $\alpha_{k}-m+1\ge \beta_{k}-1$.  We now construct two partitions $\gamma$ and $\delta$ from $\lambda$.  It is clear that $k\geq 3$. So we may define
\begin{equation}\label{con-phi-5-2-a}
  \gamma=(\beta_1+m-1,\,\ldots ,\, \beta_{k-1}+m-1 ,\,\alpha_{k+1}+1,\,\ldots,\,\alpha_{s}+1,\,2,\,1^{t-s-1})
\end{equation}
  and
\begin{equation}\label{con-phi-5-2-b}
  \delta=(\alpha_3+1-m,\,\ldots ,\, \alpha_{k}+1-m,\,\beta_{k}-1,\,\ldots ,\,\beta_{t}-1).
\end{equation}
To avoid ambiguity for the case $k=s$, we set
\begin{equation*}\label{con-phi-5-2-a1}
  \gamma=(\beta_1+m-1,\,\ldots ,\, \beta_{s-1}+m-1 ,\,2,\,1^{t-s-1}),
\end{equation*}
when $k=s$.

Using the argument in the proof of Lemma
\ref{phi-5}, it is not difficult to show that $(\gamma,\delta)_{(m+j+1)\times (j+1)}$ is an $m$-Durfee rectangle symbol. Define
 \[\phi_6(\lambda)=\left(\begin{array}{cc}
 \gamma\\[3pt]
\delta
\end{array}\right)_{(m+j+1)\times (j+1)}.\]
We claim that $\phi_6(\lambda)$ is an $m$-Durfee rectangle symbol in $P_6(-m,n)$. It is clear from \eqref{con-phi-5-2-a} and \eqref{con-phi-5-2-b} that   $\gamma_1=j+m-1$, $\delta_1=j+1$ and $\ell(\gamma)=\ell(\delta)=t-1$. It remains  to check that $|\phi_6(\lambda) |=|\lambda|$.
 Observe that
 \begin{eqnarray}
 |\gamma|+|\delta|&=&| {\alpha}|- {\alpha}_1-\alpha_2+(1-m)  (k-2)+(s-k)+2+(t-s-1)\nonumber \\[3pt]
 &&\quad +| {\beta}|+(m -1)(k-1)-(t-k+1) \nonumber\\[3pt]
&=&|\alpha|+|\beta|-{\alpha}_1-\alpha_2+m-1. \label{eqn-we-phi-6}
 \end{eqnarray}
By the definition of $Q_6(m,n)$, we have $\alpha_1=\alpha_2=j+m$. Thus from \eqref{eqn-we-phi-6},  it follows that
 \[|\gamma|+|\delta|=|\alpha|+|\beta|-(2j+m+1).
 \]
 Hence,
\begin{eqnarray*}
|\phi_6(\lambda)|&=&|\gamma|+|\delta|+(j+1)(j+m+1)\\[3pt]
&=&|\alpha|+|\beta|+j(j+m),
\end{eqnarray*}
which equals to $|\lambda|$. This proves that
$\phi_6(\lambda) \in P_6(-m,n)$.

Next we proceed to show that  $\phi_6$ is an injection. Let
\begin{equation*}
K(m,n)=\{\phi_6(\lambda)\colon \lambda\in Q_6(m,n)\}
\end{equation*}
be the set of images of $\phi_6$,
which has been shown to be a subset of $P_6(-m,n)$.
It suffices to construct a map
$\chi$ from $K(m,n)$ to $Q_6(m,n)$ such that
for any $\lambda$ in $Q_6(m,n)$, we have
\begin{equation*}\label{r-phi5-varphi}
\chi(\phi_6(\lambda))=\lambda.
\end{equation*}
To describe the map $\chi$, let
\begin{equation}\label{defi-phi-6-mu} \mu= \left(\begin{array}{cc}
 \gamma\\[3pt]
\delta
\end{array}\right)_{(m+j')\times j'}=\left(\begin{array}{ccccccc}
\gamma_1,&\gamma_2,&\ldots, &  \gamma_{t'}\\[3pt]
\delta_1,& \delta_2,&\ldots, & \delta_{t'}
\end{array}
\right)_{(m+j')\times j'}
\end{equation}
be an $m$-Durfee rectangle symbol in $K(m,n)$, that is, there is an $m$-Durfee rectangle symbol $\lambda=(\alpha,\beta)_{(m+j)\times j}$ in $Q_6(m,n)$ such that $\phi_6(\lambda)=\mu$. We claim that for $\mu \in K(m,n)$ given by \eqref{defi-phi-6-mu}, $\gamma$ has a part equal to $2$ and there exists an integer $k'$ such that
\begin{equation}\label{cond-5-k'-2}
1\leq k'\leq t'-1, \quad \gamma_{k'}-m\geq \delta_{k'}, \quad   \gamma_{k'+1}\geq 2.
\end{equation}
From the defining relation \eqref{con-phi-5-2-a} of $\phi_6$, we see that $\gamma_{s}=2$. Furthermore, by \eqref{con-phi-5-2-a} and \eqref{con-phi-5-2-b},  we have
\[2\leq k-1\leq s-1\leq t-2=t'-1,\quad \gamma_{k-1}=\beta_{k-1}+m-1,\quad \delta_{k-1}=\beta_{k}-1\]
which implies that
\[1\leq k-1\leq t'-1,\quad \gamma_{k-1}-m\ge\delta_{k-1},\quad \gamma_k\geq \gamma_s=2. \]
Hence $k-1$ satisfies the conditions in \eqref{cond-5-k'-2}. So the claim is proved.

Now we may choose $k'$ to be the maximum integer that satisfies \eqref{cond-5-k'-2}.
 From the above claim, we see that $\gamma$ has a part equal to $2$.
 Let $\gamma_{s'}=2>\gamma_{s'+1}$. The choice of $k'$ implies that  $\gamma_{k'+1}-m< \delta_{k'+1}$ when $1\leq k'\leq t'-2$.
It follows that $\delta_{k'-1}+m> \gamma_{k'+1}$ when $1\leq k'\leq t'-2$. When $k'=t'-1$, by the above claim, we see that $\gamma_{t'}\leq 2$,
which leads to $\delta_{t'-2}+m\geq \gamma_{t'}$. Combining the above two cases for $k'$, we deduce that
\begin{equation}\label{phi-6-delta-k'}
\delta_{k'-1}+m\geq \gamma_{k'+1}.
 \end{equation}
 By \eqref{cond-5-k'-2} and \eqref{phi-6-delta-k'}, we may define
 \begin{equation*}
\chi(\mu)=\left(\begin{array}{cc}
{\alpha}\\[3pt]
{\beta}
\end{array}\right)_{(m+j'-1)\times (j'-1)}.
\end{equation*}
where
\begin{equation}
\alpha=(j'+m-1,\,j'+m-1,\,\delta_1-1+m,\,\ldots, \,\delta_{k'-1}-1+m,\,\gamma_{k'+1}-1,\,\ldots,\, \gamma_{s'-1}-1)
\end{equation}
and
\begin{equation}
\beta=(\gamma_1+1-m,\,\ldots,\,\gamma_{k'}+1-m,\,\delta_{k'}+1,\, \ldots,\,\delta_{t'}+1).
\end{equation}
 It can be easily checked that $\chi(\mu)\in Q_6(m,n)$.

Finally, we are ready to verify that  $\chi(\phi_6(\lambda))=\lambda$.
By the constructions of $\phi_6(\lambda)$ and $\chi(\mu)$,   it suffices  to show that the integer $k$ appearing
 in the representation of $\phi_6(\lambda)$ is equal to the
 integer $k'$ appearing in
the representation of $\chi(\phi_6(\lambda))$ plus $1$, that is, $k'=k-1.$ This assertion can be justified by  using the same arguments as in the proof of Lemma \ref{phi-5}. For completeness, we include the following detailed proof.

Recall that  $k$ is the maximum integer determined by  $\lambda$  subject to the conditions
\begin{equation}\label{cond-k-6}
3\le k\le s, \quad  \alpha_{k}-m+1\ge \beta_{k}-1.
\end{equation}
We proceed to show that $k-1$ is the maximum integer $k'$ determined by $\phi_6(\lambda)$   such that
 \begin{equation}\label{cond-k'-6}
1\leq k'\leq t'-1, \quad \gamma_{k'}-m\geq \delta_{k'}, \quad   \gamma_{k'+1}\geq 2.
\end{equation}
From the constructions \eqref{con-phi-5-2-a} and \eqref{con-phi-5-2-b}  of $\phi_6$, it can be checked that $k-1$ satisfies the conditions in \eqref{cond-k'-6}.
So it suffices to show that $k-1$ is   the maximum integer satisfying conditions in \eqref{cond-k'-6}.  Assume to the contrary that there is an
 integer  $k \leq p \leq t'-1$  for which the conditions in \eqref{cond-k'-6}
 are satisfied, that is,
\begin{equation}\label{phi-6-gamma-p}
 \gamma_{p+1}\geq 2, \quad \gamma_{p}-m\geq \delta_{p}.
 \end{equation}
From the construction \eqref{con-phi-5-2-a} of $\gamma$, we see that $\gamma_s=2>\gamma_{s+1}$ which implies  $p+1\leq s$ under the assumption $\gamma_{p+1}\geq 2$. This yields
 \begin{equation*}\label{phi-6-k-p}
 k\leq p \leq s-1.
 \end{equation*}
 Again, in view of the constructions \eqref{con-phi-5-2-a} and  \eqref{con-phi-5-2-b} of $\phi_6$,    we find that
 \begin{equation}\label{phi-6-gam-alph}
 \gamma_{p}=\alpha_{p+1}+1,\quad \delta_{p}=\beta_{p+1}-1.
 \end{equation}
Substituting  \eqref{phi-6-gam-alph} into  \eqref{phi-6-gamma-p}, we arrive at
 \begin{equation*}
 {\alpha}_{p+1}-m+1\geq {\beta}_{p+1}-1.
 \end{equation*}
 This means that  \eqref{cond-k-6} is valid with $k$ being replaced by $p+1$. But this contradicts the maximality of $k$. So we conclude that
 $k-1$ is  the maximum integer satisfying conditions in \eqref{cond-k'-6},
 which implies that $\chi(\phi_6(\lambda))=\lambda$.
  This completes the proof. \qed

For  example, for $m=2$ and $n=60$,  let
\[\lambda=\left(\begin{array}{ccccccccc }
5,&5,&5,&5,&3,&2\\[3pt]
3,&3,&3,&3,&2,&2,&2,&2
\end{array}
\right)_{5\times 3}\]
be a $2$-Durfee rectangle symbol in $Q_6(2,60)$. It can be checked that $k=6$. Applying  $\phi_6$ to $\lambda$, we obtain
\[\mu=\phi_6(\lambda)=\left(\begin{array}{ccccccc }
4,&4,&4,&4,&3,&2,&1\\[3pt]
4,&4,&2,&1,&1,&1,&1
\end{array}
\right)_{6\times 4},\]
which  in $P_6(-2,60)$. Applying $\chi$ to $\mu$, we have $k'=5$ and $\chi(\mu)=\lambda$.

It should be noted that the injection $\phi_6$ is not valid for the case $m=0$. To be more specific, $\phi_6$ does not apply to the  Durfee  symbols $\lambda=(\alpha,\beta)_j$ in $Q_6(0,n)$ with $\beta_{s-1}=2$, where $\ell(\alpha)=s$ and $\ell(\beta)=t$. Assume that $(\gamma,\delta)_{j'}=\phi_6(\lambda)$.
 Since $\beta_{s-1}=2$, we have  $\alpha_s+1\ge 2>\beta_{s}-1$,
 which implies that $k=s$. Thus
\[
\gamma=(\beta_1-1,\,\ldots ,\, \beta_{s-1}-1 ,\,2,\,1^{t-s-1}),
\]
which is not a partition, since $\gamma_{s-1}=\beta_{s-1}-1=1$ and $\gamma_s=2$.

Combining the above injections  $\phi_i$ $(1\leq i \leq 6)$,
we are led to an injection from the set $Q(m,n)$ to the set $P(-m,n)$ for the case $m
\geq 1$.

\noindent{\it Proof of Theorem \ref{thm-equi-co} for $m\geq 1$.}
Suppose  that $m\geq 1$.  From the definitions of $Q_i(m,n)$ and $P_i(-m,n)$, we have
\[Q(m,n)=\bigcup_{i=1}^6Q_i(m,n)\]
and
\[ P(-m,n)=\bigcup_{i=1}^8P_{i}(-m,n).\]
Let $\lambda \in Q(m,n)$, define
\[\Phi(\lambda)=\begin{cases}
\phi_1(\lambda), \quad \text{if} \quad \lambda \in Q_1(m,n);\\[3pt]
\phi_2(\lambda), \quad \text{if} \quad \lambda \in Q_2(m,n);\\[3pt]
\phi_3(\lambda), \quad \text{if} \quad \lambda \in Q_3(m,n);\\[3pt]
\phi_4(\lambda), \quad \text{if} \quad \lambda \in Q_4(m,n);\\[3pt]
\phi_5(\lambda), \quad \text{if} \quad \lambda \in Q_5(m,n);\\[3pt]
\phi_6(\lambda), \quad \text{if} \quad \lambda \in Q_6(m,n).
 \end{cases}\]
Combining  Lemma \ref{phi-2} to Lemma \ref{phi-6},
 we conclude that $\Phi$ is an injection from the set $Q(m,n)$ to the set $P(-m,n)$.\qed

\section{Proof of Theorem \ref{main-equi} for $m=0$}

In this section, we  give a   proof of Theorem \ref{thm-equi-co} for   $m=0$.  In addition to the injections in Section 4,
this seemingly special case requires five more injections.

Recall that $Q(0,n)$ denotes the set of Durfee symbols $(\alpha,\beta)_j$ of $n$ such that   $\beta_1=j$ and $P(0,n)$ denotes the set of  Durfee symbols $(\gamma,\delta)_j$ of $n$ such that   $\ell(\delta)-\ell(\gamma)\leq 0.$
 From the definitions of $Q_i(0,n)$ and $P_i(0,n)$ given in Section 4, it
 can be seen that
\[Q(0,n)=\bigcup_{i=1}^6Q_i(0,n)\]
and
\[ P(0,n)=\bigcup_{i=1}^8P_{i}(0,n).\]
It is known that $Q_1(0,n)=P_1(0,n)$.
By Lemmas \ref{phi-2}, \ref{phi-3} and \ref{phi-4}, we see that the injections $ \phi_2,\phi_3, \phi_4$
can   be applied into the sets $ Q_2(0,n),\,Q_3(0,n)$ and $Q_4(0,n)$.

 As mentioned in the previous section, the injections $\phi_5$ and $\phi_6$ do not apply
  to $Q_5(0,n)$ and $Q_6(0,n)$. We need to construct an injection from  $Q_5(0,n)\cup Q_6(0,n)$ to $P_5(0,n)\cup P_6(0,n)\cup P_7(0,n)\cup P_8(0,n)$. To this end,  we shall divide the set $ Q_5(0,n)\cup Q_6(0,n)$ into the following five disjoint subsets $\bar{Q}_1(0,n)$,
 $\bar{Q}_2(0,n)$, $\bar{Q}_3(0,n)$, $\bar{Q}_4(0,n)$ and $\bar{Q}_5(0,n)$:
\begin{itemize}

\item[(1)] $\bar{Q}_1(0,n)$ is the set of Durfee symbols $(\alpha,\beta)_{j} \in Q_5(0,n)$ with $s(\beta)\geq 3$;

\item[(2)] $\bar{Q}_2(0,n)$ is the set of Durfee symbols $(\alpha,\beta)_{j} \in  Q_6(0,n)$ with $s(\beta)\geq 3$;

\item[(3)] $\bar{Q}_3(0,n)$ is the set of Durfee symbols $(\alpha,\beta)_{j} \in  Q_5(0,n)\cup Q_6(0,n)$ with $s(\alpha)=1$ and $s(\beta)= 2$;

\item[(4)] $\bar{Q}_4(0,n)$ is the set of  Durfee symbols $(\alpha,\beta)_{j} \in Q_5(0,n)\cup Q_6(0,n)$ with $s(\alpha)\geq 2$, $\beta_1=\beta_2$ and $s(\beta)= 2$;

\item[(5)] $\bar{Q}_5(0,n)$ is the set of  Durfee symbols $(\alpha,\beta)_{j} \in Q_5(0,n)\cup Q_6(0,n)$ with $s(\alpha)\geq 2$, $\beta_1>\beta_2$ and $s(\beta)= 2$.

\end{itemize}

On the other hand, we  divide the set  $P_5(0,n)\cup P_6(0,n)$ into   three disjoint subsets $\bar{P}_1(0,n)$, $\bar{P}_2(0,n)$ and $\bar{P}_3(0,n)$:
\begin{itemize}
\item[(1)] $\bar{P}_1(0,n)$ is the set of  Durfee symbols $(\gamma,\delta)_{j'} \in P_5(0,n)$ with $s(\delta)\geq 2$;

\item[(2)] $\bar{P}_2(0,n)$ is the set of  Durfee symbols $(\gamma,\delta)_{j'} \in P_6(0,n)$ with $s(\delta)\geq 2$;

\item[(3)] $\bar{P}_3(0,n)$ is the set of  Durfee symbols $(\gamma,\delta)_{j'} \in P_5(0,n)\cup P_6(0,n)$ with $s(\delta)=1$.
\end{itemize}

In the following lemmas, we shall show that
 there exist an injection $\psi_1$   from  $\bar{Q}_1(0,n)$ to  $\bar{P}_1(0,n)$,  an injection $\psi_2$   from    $\bar{Q}_2(0,n)$ to  $\bar{P}_2(0,n)$, an injection  $\psi_3$ from   $\bar{Q}_3(0,n)$ to $\bar{P}_3(0,n)$,  an injection $\psi_4$ from   $\bar{Q}_4(0,n)$ to  $P_7(0,n)$ and an injection $\psi_5$ from  $\bar{Q}_5(0,n)$ to $P_8(0,n)$.  Then the injection $\Phi$ for $m=0$
 can be composed from injections  $\phi_i$ $(1\leq i\leq 4)$ and  injections $\psi_i$  $(1\leq i\leq 5)$.

 \begin{lem}\label{psi-1}
There exists an injection $\psi_1$   from   the set $\bar{Q}_1(0,n)$ to the set $\bar{P}_1(0,n)$.
\end{lem}
\pf Let
 \[\lambda=\left(\begin{array}{c }
\alpha\\[3pt]
\beta
\end{array}
\right)_{j}=\left(\begin{array}{cccc }
\alpha_1,&\alpha_2,&\ldots, & \alpha_s\\[3pt]
\beta_1,& \beta_2,&\ldots, & \beta_t
\end{array}
\right)_{j}\]
be a Durfee  symbol in $\bar{Q}_1(0,n)$. By definition, we have $j=\beta_1\geq \beta_t\geq 3$, $j=\alpha_1=\alpha_2>\alpha_3$ and $t-s\geq 1$.
Consequently, we have $\alpha_2+2=j+2 > \beta_3-1$.
 Hence there exists the maximum number $k$ such that $1\leq k\leq t-1$ and $\alpha_k+2\ge \beta_{k+1}-1$. So we may define
  \begin{equation}\label{psi-0-lambda}
  \psi_1(\lambda)=\left(\begin{array}{cc}
 \gamma\\[3pt]
\delta
\end{array}\right)_{j+1},
\end{equation}
where
\[
\gamma=(\beta_2-2,\,\ldots,\, \beta_k-2,\,\alpha_{k+1}+1,\,\ldots,\,\alpha_{t}+1)
\]
and
\[
\delta=(\alpha_2+1,\,\alpha_3+2,\,\ldots, \, \alpha_k+2,\,\beta_{k+1}-1,\,\ldots,\,\beta_{t}-1).
\]
Using the same arguments as in the proof of Lemma \ref{phi-5}, we
deduce that $\psi_1(\lambda)$ is a Durfee symbol in $\bar{P}_1(0,n)$ and the construction of $\psi_1$ is
reversible. Hence  $\psi_1$ is an injection from $\bar{Q}_1(0,n)$ to  $\bar{P}_1(0,n)$. This completes the proof. \qed

  \begin{lem}\label{psi-2}
There exists an injection $\psi_2$   from   the set $\bar{Q}_2(0,n)$ to the set $\bar{P}_2(0,n)$.
\end{lem}
\pf Let
 \[\lambda=\left(\begin{array}{c }
\alpha\\[3pt]
\beta
\end{array}
\right)_{j}=\left(\begin{array}{cccc }
\alpha_1,&\alpha_2,&\ldots, & \alpha_s\\[3pt]
\beta_1,& \beta_2,&\ldots, & \beta_t
\end{array}
\right)_{j}\]
be a Durfee  symbol in $\bar{Q}_2(0,n)$. By definition, we have $j=\beta_1\geq \beta_t\geq 3$, $j=\alpha_1=\alpha_2=\alpha_3$ and $t-s\geq 1$.
Thus, $\alpha_3+1=j+1>\beta_3-1$. So we may choose
  the maximum integer $k$ such that $  k\le s$ and $\alpha_{k}+1\ge \beta_{k}-1$. Define
  \begin{equation}\label{psi-0-lambda}
  \psi_2(\lambda)=\left(\begin{array}{cc}
 \gamma\\[3pt]
\delta
\end{array}\right)_{j+1},
\end{equation}
where
\[
  \gamma=(\beta_1-1,\,\ldots ,\, \beta_{k-1}-1 ,\,\alpha_{k+1}+1,\,\ldots,\,\alpha_{s}+1,\,2,\,1^{t-s-1})
\]
  and
\[
  \delta=(\alpha_3+1,\,\ldots ,\, \alpha_{k}+1,\,\beta_{k}-1,\,\ldots ,\,\beta_{t}-1).
\]
It can be checked that $\psi_2(\lambda)$ is a Durfee symbol in $\bar{P}_2(0,n)$ and the construction of $\psi_2$ is
reversible by the same argument  as in the proof of Lemma \ref{phi-6}.
So we deduce that $\psi_2$ is an injection. This completes the proof.  \qed

\begin{lem} \label{psi-3}
There is an injection $\psi_3$ from the set $\bar{Q}_3(0,n)$ to the set $\bar{P}_3(0,n)$.
\end{lem}
\pf
Let
\[\lambda=\left(\begin{array}{cc}
\alpha\\[3pt]
\beta
\end{array}\right)_{j}=\left(\begin{array}{cccc }
\alpha_1,&\alpha_2,&\ldots ,& \alpha_s\\[3pt]
\beta_1,& \beta_2,&\ldots ,& \beta_t
\end{array}
\right)_{j} \]
be a Durfee symbol in $\bar{Q}_3(0,n)$. By definition $\alpha_1=\alpha_2=j$,  $\alpha_s=1$, $\beta_1=j$, $\beta_t=2$ and $t-s\ge 1$.
 So we have $\beta_2\leq j$ and $\alpha_2=j$. This enables us to define
\begin{equation}\label{psi-1-con}
\psi_3(\lambda)=\left(\begin{array}{cc}
\gamma\\[3pt]
\delta
\end{array}\right)_{j'}=\left(\begin{array}{cccccc cc}
  \beta_2-1,&\ldots, &\beta_t-1\\[3pt]
  \alpha_2+1,&\ldots,&\alpha_{s-1}+1,&1^{t-s+1}
\end{array}
\right)_{j+1}.
\end{equation}
 Note that $\gamma_1=\beta_2-1\le j-1=j'-2$, $\delta_1=\alpha_2+1=j+1=j'$ and $\ell(\gamma)=\ell(\delta)$. Since $t-s\geq 1$, we see that $s(\delta)=1$. Since $\alpha_s=1$, it is easy to check that $|\psi_3(\lambda)|=|\lambda|$.  So  $\psi_3(\lambda)$ is in $ \bar{P}_3(0,n)$.

We now proceed to show that $\psi_3$ is an injection.
Let
\begin{equation*}
L(m,n)=\{\psi_3(\lambda)\colon \lambda\in \bar{Q}_3(0,n)\}
\end{equation*}
be the set of images of $\psi_3$,
which has been shown to be a subset of $\bar{P}_3(0,n)$.
It suffices to construct a map
$\vartheta$ from $L(m,n)$ to $\bar{Q}_3(0,n)$ such that
for any $\lambda$ in $\bar{Q}_3(0,n)$, we have
\begin{equation}\label{r-psi-1}
\vartheta(\psi_3(\lambda))=\lambda.
\end{equation}

 Let
 \[\mu= \left(\begin{array}{cc}
 \gamma\\[3pt]
\delta
\end{array}\right)_{j'}=\left(\begin{array}{ccccccc}
\gamma_1,&\gamma_2,&\ldots, &  \gamma_{t'}\\[3pt]
\delta_1,& \delta_2,&\ldots, & \delta_{t'}
\end{array}
\right)_{j'}\]
be a Durfee symbol in $L(m,n)$. We claim that $\gamma_{t'}=1$ and $\delta_{t'-1} =1.$ By the definition of $L(m,n)$,  there exists $\lambda$ in $\bar{Q}_3(0,n)$ such that $\psi_3(\lambda)=\mu$. Since $t-s+1\ge 2$ and $\beta_t=2$, from
the definition \eqref{psi-1-con} of $\psi_3(\lambda)$, we get
\begin{equation}\label{psi-1-a}
\gamma_{t'}=\beta_{t}-1=1, \quad \delta_{t'-1}=1.
\end{equation}
So the claim holds.

We next define the map $\vartheta$. Let ${h'}$ be the largest index such that $\delta_{h'}>1$. By the above claim,  we have $\delta_{t'-1}=1$, and so $h'\leq t'-2$.
Define
\[\vartheta(\mu)=\left(
\begin{array}{ccccccccc}
j'-1,&\delta_1-1,&\ldots,&\delta_{h'}-1,&1\\[3pt]
j'-1,&\gamma_1+1,&\ldots,&\gamma_{t'}+1&
\end{array}
\right)_{j'-1}.\]
It is easy to verify that $\vartheta(\mu)\in \bar{Q}_3(0,n)$ and  $\vartheta(\psi_3(\lambda))=\lambda$ for $\lambda \in \bar{Q}_3(0,n)$. Therefore,   $\psi_3$ is an injection from $\bar{Q}_3(0,n)$ to $\bar{P}_3(0,n)$. This completes the proof. \qed

For example, for $n=35$, let
\[\lambda=\left(\begin{array}{ccccccc }
3,&3,&2,&2,&1\\[3pt]
3,&3,&3,&2,&2,&2
\end{array}
\right)_3,\]
which is a Durfee symbol in $\bar{Q}_3(0,35)$.
Applying the injection $\psi_3$ to $\lambda$, we get
\[\psi_3(\lambda)=\left(\begin{array}{ccccccc }
2,&2,&1,&1,&1\\[3pt]
4,&3,&3,&1,&1
\end{array}
\right)_4,\]
which is in  $\bar{P}_3(0,35)$. Applying $\vartheta$ to $\mu$, we recover $\lambda$, that is,  $\vartheta(\mu)=\lambda$.

\begin{lem}\label{psi-4}
There is a bijection $\psi_4$ between the set $\bar{Q}_4(0,n)$ and the set $P_7(0,n)$.
\end{lem}
\pf Let
\[\lambda=\left(\begin{array}{cc}
\alpha\\[3pt]
\beta
\end{array}\right)_j=\left(\begin{array}{cccc }
\alpha_1,&\alpha_2,&\ldots ,& \alpha_s\\[3pt]
\beta_1,& \beta_2,&\ldots ,& \beta_t
\end{array}
\right)_j\]
be a Durfee symbol  in $\bar{Q}_4(0,n)$. By definition, $\alpha_1=\alpha_2=j$, $\alpha_s\ge 2$, $\beta_1=\beta_2=j$,  $\beta_t=2$ and $t-s\ge 1$. Thus,  $\alpha_2=j>\beta_3-1$ and $\beta_2=j\geq \alpha_3$. So we may define
 \[\psi_4(\lambda)=\left(\begin{array}{cc}
\gamma\\[3pt]
\delta
\end{array}\right)_{j'}=\left(\begin{array}{cccccc cc}
 \alpha_2,&\beta_3-1,&\ldots,&\beta_{t-1}-1\\[3pt]
 \beta_2+1,& \alpha_3+1,&\ldots,&\alpha_{s}+1, & 1^{t-s-1}
\end{array}
\right)_{j+1}.\]
 Note that $\delta_{s-1}=\alpha_s+1\ge 3$ and $\delta_i=1$ for $s\le i\le t-2$.
 It is clear that $\delta$ has no parts equal to $2$.
   Since $\beta_t=2$, we find that $|\psi_4(\lambda)|=|\lambda|$.
   Moreover, we have  $\ell(\gamma)=\ell(\delta)=t-2$, $\delta_1=\beta_2+1=j+1=j'$ and $\gamma_1=\alpha_2=j=j'-1>\beta_3-1=\gamma_2$. So $\psi_4(\lambda)$ is in $P_7(0,n)$.

   It is not difficult to verify that
   $\psi_4$ is invertible. Thus $\psi_4$ is a  bijection.
    \qed

For example, for $n=40$, consider the following Durfee symbol in $\bar{Q}_4(0,40)$:
\[\lambda=\left(\begin{array}{ccccccc }
3,&3,&3,&2,&2\\[3pt]
3,&3,&3,&3,&2,&2,&2
\end{array}
\right)_{3}.\]
Applying the bijection $\psi_4$, we get
\[\psi_4(\lambda)=\left(\begin{array}{ccccccc }
3,&2,&2,&1,&1\\[3pt]
4,&4,&3,&3,&1
\end{array}
\right)_{4},\]
which is  in $P_7(0,40)$.

\begin{lem}\label{psi-5}
There is an injection $\psi_5$ from the set $\bar{Q}_5(0,n)$ to the set $P_8(0,n)$.
\end{lem}

\pf  Let
\[\lambda=\left(\begin{array}{cc}
\alpha\\[3pt]
\beta
\end{array}\right)_j=\left(\begin{array}{cccc }
\alpha_1,&\alpha_2,&\ldots ,& \alpha_s\\[3pt]
\beta_1,& \beta_2,&\ldots ,& \beta_t
\end{array}
\right)_j\]
be a Durfee symbol  in $\bar{Q}_5(0,n)$. By definition, $\alpha_1=\alpha_2=j$, $\alpha_s\ge 2$, $j=\beta_1>\beta_2$,  $\beta_t=2$ and $t-s\ge 1$.
 In this case, we have $j\geq 3$. By definition, $j\geq 2$.  If $j=2$, then $\beta_1=\beta_2=\cdots=\beta_t=2$, which is a contradiction
  to $\beta_1>\beta_2$.

We next define the map $\psi_5$. Since $\alpha_1=\alpha_2=j$,
 we may choose $k$ to be the maximum number such that $\alpha_k=j$. Clearly,  $k\ge 2$. Since $\beta_t=2$, we may choose  $h$ to be the minimum number such that $\beta_h=2$. Since $\beta_1=j>2$, we get $2\leq h\leq t$.

  By the choice of $k$ and $j=\beta_1>\beta_2$, we see that $\alpha_k>\beta_2$ and $\beta_1>\alpha_{k+1}$. On the other hand, by the choice of $h$, we see that $\beta_{h-1}>\beta_h$. So we may define
 \begin{align}\label{psi-5-co}
\psi_5(\lambda)&= \left(\begin{array}{cc}
 \gamma\\[3pt]
\delta
\end{array}\right)_{j'} \nonumber \\[3pt]
&=\left(\begin{array}{ccccccccc cc}
\alpha_1-1,&\ldots,&\alpha_k-1,&\beta_2-1,&\ldots,&\beta_{h-1}-1,&\beta_h, &1^{t-h}\\[3pt]
\beta_1,&\alpha_{k+1}+1,& \ldots,&\alpha_s+1,&1^{2k-2+t-s}
\end{array}
\right)_{j}.
\end{align}
 Since $\delta_{s-k+1}=\alpha_{s}+1\ge 3$ and $\delta_i=1$ for $s-k+2\le i\le t-1+k$, we deduce that $\delta$ has no parts equal to  $2$. Furthermore,
 it is easily checked that  $\ell(\gamma)=\ell(\delta)=t+k-1$, $\delta_1=j'$ and $\gamma_1=\gamma_2=j'-1$ and $|\psi_5(\lambda)|=|\lambda|$. So $\psi_5(\lambda)$ is in $P_8(0,n)$.

To prove that $\psi_5$ is an injection, let
\begin{equation*}
R(m,n)=\{\psi_5(\lambda)\colon \lambda\in \bar{Q}_5(0,n)\}
\end{equation*}
be the set of images of $\psi_5$,
which has been shown to be a subset of $P_8(0,n)$.
We shall construct a map
$\theta$ from $R(m,n)$ to $\bar{Q}_5(0,n)$ such that
for any $\lambda$ in $\bar{Q}_5(0,n)$, we have
\begin{equation*}\label{r-psi-3}
\theta(\psi_5(\lambda))=\lambda.
\end{equation*}

 Let
 \[\mu= \left(\begin{array}{cc}
 \gamma\\[3pt]
\delta
\end{array}\right)_{j'}=\left(\begin{array}{ccccccc}
\gamma_1,&\gamma_2,&\ldots, &  \gamma_{t'}\\[3pt]
\delta_1,& \delta_2,&\ldots, & \delta_{t'}
\end{array}
\right)_{j'}\]
be a Durfee symbol in $R(m,n)$. By the definition of $R(m,n)$, we see that there exists $\lambda$ in $\bar{Q}_5(0,n)$ such that $\psi_5(\lambda)=\mu$. Let $k'$ denote the number of parts in $\gamma$ equal to $j'-1$ and $n_1(\delta)$ denote the number of parts $1$ in $\delta$. We claim that if $j'\ge 4$, then  $n_1(\delta)\geq 2k'-1$; and if $j'=3$, then $k'\geq 3$ and $n_1(\delta)\geq 2k'-3$.

 From the construction \eqref{psi-5-co}  of $\psi_5$, we find that $j'=j$ and $n_1(\delta)=2k-2+t-s$. Since $t-s\geq 1$, we get  $n_1(\delta)\geq 2k-1$.
 Moreover, since $k\geq 2$, it suffices to show that $k'=k$ if $j\geq 4$ and $k'=k+1$ if $j=3$.
   From the construction \eqref{psi-5-co}  of $\psi_5$, we have $\gamma_i=\alpha_i-1$ for $1\leq i\leq k$. Since $\alpha_i=j$  for  $1\leq i\leq k$, it follows that $\gamma_i=j-1$ for $1\leq i \leq k$.

 It remains to show that $\gamma_{k+1}<j-1$ if $j\geq 4$ and $\gamma_{k+1}=j-1>\gamma_{k+2}$ if $j=3$.   By \eqref{psi-5-co}, we have either $\gamma_{k+1}=\beta_2-1$ or $\gamma_{k+1}=\beta_h$. If $j\geq 4$, in either case, we find that $\gamma_{k+1}<j-1$ since $\beta_2<j$ and $\beta_h=2$.  If $j=3$, then $\beta_1=3$ and $\beta_2=2$, so $h=2$ where we recall that $h$ be the minimum integer such that $\beta_h=2$. This implies that $\gamma_{k+1}=\beta_2=2=j-1$. Since $\gamma_{k+2}\leq 1$, we get  $\gamma_{k+2}< j-1$. Thus, we arrive at the conclusion that $k'=k+1$. So the claim is verified.

 From the construction \eqref{psi-5-co}  of $\psi_5$, it can be seen that $\gamma_{k+h-1}=\beta_h=2$. We may choose $h'$ to be the maximum integer such that $\gamma_{h'}=2$.  Assumed that   $\gamma$ has $k'$ parts equal to $j'-1$. From the above claim, we see that if $j'\geq 4$, then $k'\geq 2$. If $j'=3$, then $k'\geq 3$.  We consider the following two cases:

\noindent{Case 1: $j'\ge 4$.}  Define
\begin{align*}
\theta(\mu)
=\left(\begin{array}{ccccccccc }
\gamma_1+1,&\ldots,&\gamma_{k'}+1,&\delta_2-1,&\ldots,&\delta_{t'}-1\\[3pt]
\delta_1,&\gamma_{k'+1}+1,&\ldots,&\gamma_{h'-1}+1,& \gamma_{h'},&\gamma_{h'+1}+1,&\ldots,&\gamma_{t'}+1
\end{array}
\right)_{j'}.\end{align*}
By the above claim, we have $n_1(\delta)\geq 2k'-1$. Now it is
easy to check that $\theta(\mu)\in \bar{Q}_5(0,n)$.

\noindent{Case 2:} $j'=3$. By the definition of $k'$ and $h'$, we have $k'=h'$. Let $r'=n_1(\delta)$, and  define
 \[\theta(\mu)
=\left(\begin{array}{ccccccccc cc}
3^{k'-1},&2^{t'-r'-1}\\[3pt]
3,&2^{t'-k'+1}
\end{array}
\right)_{3}.
\]
From the above claim, we have $r'\geq 2k'-3$ and $k'\geq 3$.
Now it can be easily checked that  $\theta(\mu) \in \bar{Q}_5(0,n)$.

Finally, from the constructions of $\psi_5$ and $\theta$ and the above claim, it is straightforward to verify that $\theta(\psi_5(\lambda))=\lambda$ for any $\lambda\in \bar{Q}_5(0,n)$. This completes the proof. \qed

For example, for $n=51$, consider the following Durfee symbol in $\bar{Q}_5(0,51)$:
\[\lambda=\left(\begin{array}{ccccccc }
4,&4,&4,&2,&2\\[3pt]
4,&3,&3,&3,&2,&2,&2
\end{array}
\right)_{4}.\]
Applying the injection $\psi_5$,  we see that $k=3$ and $h=4$, and we obtain
\[\mu=\psi_5(\lambda)=\left(\begin{array}{ccccccccccccc }
3,&3,&3,&2,&2,&2,&2,&1,&1\\[3pt]
4,&3,&3,&1,&1,&1,&1,&1,&1
\end{array}
\right)_{4},\]
which is  in $P_8(0,51)$. Applying $\theta$ to $\mu$, we find that $k'=3$, $h'=4$ and $\theta(\mu)=\lambda$.

We are now ready to complete the proof of Theorem \ref{thm-equi-co} for the
case $m=0$.

\noindent{\it Proof of Theorem \ref{thm-equi-co} for $m=0$.}
From the definitions of $Q_i(0,n)$ ($1\leq i\leq 4$) and $\bar{Q}_i(0,n)$ ($1\leq i\leq 5$),  we have
\begin{eqnarray*}
Q(0,n)&=&Q_1(0,n)\cup Q_2(0,n)\cup Q_3(0,n)\cup Q_4(0,n)\cup \bar{Q}_1(0,n)\cup\bar{Q}_2(0,n) \\[3pt]
&&\quad \cup \,\bar{Q}_3(0,n)\cup \bar{Q}_4(0,n)\cup \bar{Q}_5(0,n).
\end{eqnarray*}
\underline{}By the definitions of $P_i(0,n)$ ($1\leq i\leq 8$)  and $\bar{P}_i(0,n)$ ($1\leq i\leq 3$),  we have
\begin{eqnarray*}
 P(0,n)&=&P_1(0,n)\cup P_2(0,n)\cup P_3(0,n)\cup P_4(0,n) \cup \bar{P}_1(0,n)\cup\bar{P}_2(0,n)\\[3pt]
 &&\quad \cup\,\bar{P}_3(0,n)\cup P_7(0,n)\cup P_8(0,n).
 \end{eqnarray*}
Let $\lambda \in Q(0,n)$, define
\[\Phi(\lambda)=\begin{cases}
\phi_1(\lambda), \quad \text{if} \quad \lambda \in Q_1(0,n);\\[3pt]
\phi_2(\lambda), \quad \text{if} \quad \lambda \in Q_2(0,n);\\[3pt]
\phi_3(\lambda), \quad \text{if} \quad \lambda \in Q_3(0,n);\\[3pt]
\phi_4(\lambda), \quad \text{if} \quad \lambda \in Q_4(0,n);\\[3pt]
\psi_1(\lambda), \quad \text{if} \quad \lambda \in \bar{Q}_1(0,n);\\[3pt]
\psi_2(\lambda), \quad \text{if} \quad \lambda \in \bar{Q}_2(0,n);\\[3pt]
\psi_3(\lambda), \quad \text{if} \quad \lambda \in \bar{Q}_3(0,n);\\[3pt]
\psi_4(\lambda), \quad \text{if} \quad \lambda \in \bar{Q}_4(0,n);\\[3pt]
\psi_5(\lambda), \quad \text{if} \quad \lambda \in \bar{Q}_5(0,n).
 \end{cases}\]
 From   Lemmas \ref{phi-2} to \ref{phi-4} and Lemmas \ref{psi-1} to   \ref{psi-5}, it immediately follows that $\Phi$ is an injection from   $Q(0,n)$ to   $P(0,n)$. This completes the proof. \qed

By Theorem \ref{equi-con}, we confirm Conjecture \ref{conj} of Andrews, Dyson and Rhoades, or equivalently,  their Conjecture \ref{conj-o} on the
spt-crank.

By a closer examination of the injections in the proof of
Theorem \ref{thm-equi-co}, we can characterize the numbers
$n$ and $m$ for which $N_{\le m}(n)=M_{\le m}(n)$. Here we omit the detailed analysis.

\section{Connection to  Theorem \ref{main-c}}

In this section, we establish a connection between
 Conjecture \ref{conj} and  Theorem \ref{main-c} of Andrews, Chan and Kim.
 More precisely, we relate the positive rank (crank) moments $\overline{N}_k(n)$ ($\overline{M}_k(n)$) to the functions $N_{\leq m}(n) $ ($M_{\leq m}(n)$) defined by Andrews, Dyson and Rhoades. Based on this connection, it can be
 seen that Theorem \ref{main-c} of Andrews, Chan and Kim on the    positive rank and crank moments  can be deduced from Conjecture \ref{conj}.
 This leads to an alternative proof of the theorem of Andrews, Chan and Kim.

\begin{thm}\label{con-ra-func} For $k\geq 1$ and $n\geq 1$, we have
\begin{eqnarray}\label{las-r-r}
 \overline{N}_k(n)&=&\frac{1}{2}\sum_{m=1}^{+\infty}(m^k-(m-1)^k)\left( {p(n)-N_{\leq m-1}(n)} \right),\\[3pt]
\overline{M}_k(n)&=&\frac{1}{2}\sum_{m=1}^{+\infty}(m^k-(m-1)^k)\left( {p(n)-M_{\leq m-1}(n)} \right). \label{las-r-c}
\end{eqnarray}
\end{thm}
\pf We only give a proof of \eqref{las-r-r} since  \eqref{las-r-c}
can be justified in the same vain. Recall that
\begin{equation}\label{expre-n}
\overline{N}_k(n)=\sum_{j=1}^{+\infty}j^kN(j,n).
\end{equation}
Express \eqref{expre-n} in the following form:
\begin{align*}
\overline{N}_k(n)&=\sum_{j=1}^{+\infty} N(j,n)\left(\sum_{m=1}^j  m^k- \sum_{m=1}^j(m-1)^k\right)\\[3pt]
&=\sum_{j=1}^{+\infty}\sum_{m=1}^j (m^k-(m-1)^k) N(j,n).
\end{align*}
Changing the order of  summations, we find that
 \begin{equation}\label{last-4-a}
\overline{N}_k(n)=\sum_{m=1}^{+\infty}(m^k-(m-1)^k)\sum_{j=m}^{+\infty}N(j,n).
\end{equation}
Writing the second sum in \eqref{last-4-a} as
\begin{equation}\label{last-5}
\sum_{j=m}^{+\infty}N(j,n)=\sum_{j=-\infty}^{+\infty}N(j,n)
-\sum_{j=-\infty}^{m-1}N(j,n),
\end{equation}
and substituting the relations
\[\sum_{r=-\infty}^{\infty}N(r,n)=p(n)
\] and
\[\sum_{j=-\infty}^{m-1}N(j,n)=p(-m+1,n)\]
as given by  \eqref{sum} and \eqref{temp-2}
into \eqref{last-5}, we deduce that
\begin{eqnarray}\label{last-2-a}
\sum_{j=m}^{+\infty}N(j,n)&=&p(n)-p(-m+1,n).
\end{eqnarray}
Replacing   $m$ by $m-1$ in \eqref{Nleq} yields
\begin{equation}\label{last-1-a}
p(-m+1,n)=\frac{p(n)+N_{\leq m-1}(n)}{2}.
\end{equation}
Substituting \eqref{last-1-a} into \eqref{last-2-a}, we obtain
\begin{equation}\label{last-3-a}
\sum_{j=m}^{+\infty}N(j,n)=\frac{p(n)-N_{\leq m-1}(n)}{2}.
\end{equation}
Combining   \eqref{last-4-a} and \eqref{last-3-a},
we arrive at   relation \eqref{las-r-r}. This completes the proof.  \qed

In view of  Theorem \ref{con-ra-func}, it can be seen that
Theorem \ref{main-c} follows from Conjecture \ref{conj}.

\noindent{\it Proof of Theorem \ref{main-c}. }
Subtracting \eqref{las-r-r}  from  \eqref{las-r-c} in Theorem \ref{con-ra-func}, we obtain
\begin{align}\label{ineq-3}
\overline{M}_k(n)-\overline{N}_k(n)&=\frac{1}{2}\sum_{m=1}^{+\infty}(m^k-(m-1)^k)(N_{\leq m-1}(n)-M_{\leq m-1}(n)).
\end{align}
From the definitions of the rank and the crank, we have for $m\geq n+1$,
\begin{eqnarray*}
N_{\leq m-1}(n)&=&p(n),\\[3pt]
M_{\leq m-1}(n)&=&p(n).
\end{eqnarray*}
It follows that for  $m\geq n+1$
\begin{equation}\label{thm17-n+1}
N_{\leq m-1}(n)-M_{\leq m-1}(n)=0.
\end{equation}
For $m=n$, from the definitions of the rank and the crank, we find that
\begin{eqnarray*}
N_{\leq n-1}(n)&=&p(n),\\[3pt]
M_{\leq n-1}(n)&=&p(n)-2.
\end{eqnarray*}
Consequently,
\begin{equation}\label{thm17-n}
N_{\leq n-1}(n)-M_{\leq n-1}(n)=2.
\end{equation}
 Substituting
 \eqref{thm17-n+1} and \eqref{thm17-n} into  \eqref{ineq-3}, we obtain
\begin{align}
\overline{M}_k(n)-\overline{N}_k(n)&=\frac{1}{2}\sum_{m=1}^{n-1}(m^k-(m-1)^k)(N_{\leq m-1}(n)-M_{\leq m-1}(n))\nonumber\\[3pt]
&\qquad \quad + n^k-(n-1)^k. \label{diff-mk-nk}
\end{align}
Since  $m^k-(m-1)^k>0$ for $m\geq 1$ and $k\geq 1$,
  by Conjecture \ref{conj}, that is,  $N_{\leq m-1}(n)-M_{\leq m-1}(n)\geq 0$,
  we reach the assertion that $\overline{M}_k(n)-\overline{N}_k(n)>0$ for $n\geq 1$ and $k\geq 1$.
  This completes the proof. \qed

\vskip 0.5cm

\noindent{\bf Acknowledgments.} This work was supported by the 973 Project, the PCSIRT Project of the Ministry of Education and the National Science Foundation of China.

\end{document}